%% file: paper.tex
\documentclass[12pt]{article}

\usepackage{amsfonts, amssymb, amsmath, amsthm, epsf}

\makeatletter
\@addtoreset{figure}{section}

\newtheorem{theorem}{Theorem}
\@addtoreset{theorem}{section}

\newtheorem{lemma}{Lemma}
\@addtoreset{lemma}{section}

\@addtoreset{definition}{section}

\newtheorem{remark}{Remark}
\@addtoreset{remark}{section}

\newtheorem{example}{Example}
\@addtoreset{example}{section}

\newtheorem{condition}{Condition}
\@addtoreset{condition}{section}

\@addtoreset{corollary}{section}

\@addtoreset{equation}{section}

\renewcommand{\Im}{{\rm Im\,}}

\newcommand{\ind}{{\rm ind\,}}
\newcommand{\supp}{{\rm supp\,}}

\renewcommand{\ker}{{\rm ker\,}}
\renewcommand{\dim}{{\rm dim\,}}
\newcommand{\codim}{{\rm codim\,}}
\newcommand{\arctg}{{\rm arctan\,}}
\newcommand{\Orb}{{\rm Orb}}

\makeatother

\title{Nonlocal elliptic problems with nonlinear
argument transformations near the points of conjugation}
\author{Pavel~Gurevich\thanks{This research was supported by
Russian Foundation for Basic Research (grant No~03-01-06523),
Russian Ministry for Education (grant No~E02-1.0-131), and INTAS
(grant~YSF~2002-008).}}

\date{}

\begin{document}
\maketitle

\begin{abstract}
We consider elliptic equations of order $2m$ in a domain
$G\subset\mathbb R^n$ with nonlocal conditions that connect the
values of the unknown function and its derivatives on
$(n-1)$-dimensional submanifolds $\Upsilon_i$ (where
$\bigcup_i\Upsilon_i=\partial G$)  with the values on
$\omega_{is}(\overline\Upsilon_i)\subset\overline G$. Nonlocal
elliptic problems in dihedral angles arise as model problems near
the conjugation points
$g\in\overline\Upsilon_i\cap\Upsilon_j\ne\varnothing$, $i\ne j$.
We study the case where the transformations $\omega_{is}$
correspond to nonlinear transformations in the model problems. It
is proved that the operator of the problem remains Fredholm and
its index does not change as we pass from linear argument
transformations to nonlinear ones.
\end{abstract}

\tableofcontents

\input introd.tex
\input sect1.tex
\input sect2.tex
\input sect3.tex
\input sect4.tex
\input sect5.tex
\input sect6.tex

\input bibl.tex
\end{document}

%% file: introd.tex
\section*{Introduction}
\addcontentsline{toc}{section}{\protect\numberline{}Introduction}

The first mathematicians who studied ordinary differential
equations with nonlocal conditions were
Sommerfeld~\cite{Sommerfeld}, Tamarkin~\cite{Tamarkin},
Picone~\cite{Picone}. In~1932, Carleman~\cite{Carleman} considered
the problem of finding a holomorphic function in a bounded domain
$G$, satisfying the following condition: the value of the unknown
function at each point $x$ of the boundary is connected with the
value at $\omega(x)$, where $\omega\big(\omega(x)\big)=x$,
$\omega(\partial G)=\partial G$. Such a statement of the problem
originated further investigations of nonlocal elliptic problems
with the shifts mapping the boundary onto itself. In~1969,
Bitsadze and Samarskii~\cite{BitzSam} considered essentially
different type of nonlocal problems. They studied the Laplace
equation in a bounded domain $G$ with the boundary-value condition
connecting the values of the unknown function on a manifold
$\Upsilon_1\subset\partial G$ with the values on some manifold
inside $G$; on the set $\partial G\setminus \Upsilon_1$ the
Dirichlet condition was imposed. In a general case, such a problem
was formulated as an unsolved one.

\begin{figure}[ht]
{ \hfill\epsfbox{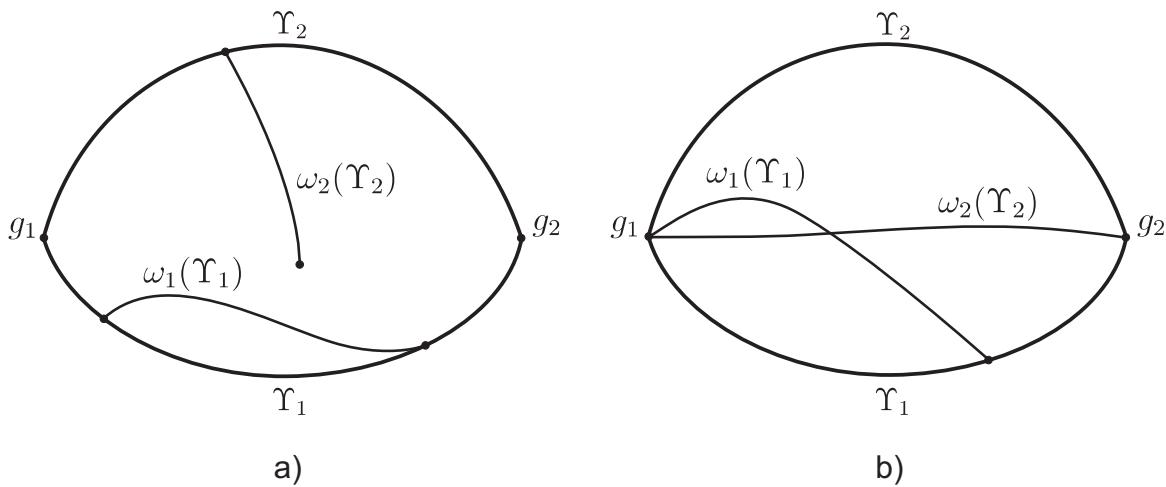}\hfill\ } \caption{The domain $G$ with
the boundary $\partial G=\bar\Upsilon_1\cup\bar\Upsilon_2$ for
$n=2$.}
   \label{figDomG}
\end{figure}

The most difficult situation in the theory of nonlocal problems is
that where the support of nonlocal terms intersects with the
boundary of domain. We consider the following example. Let
$G\subset\mathbb R^n$ ($n\ge2$) be a bounded domain with the
boundary $\partial G=\Upsilon_1\cup\Upsilon_2\cup\mathcal K_1$,
where $\Upsilon_i$ are smooth open (in the topology of $\partial
G$) $(n-1)$-dimensional $C^\infty$-manifolds, $\mathcal
K_1=\bar\Upsilon_1\cap\bar\Upsilon_2$ is an $(n-2)$-dimensional
connected $C^\infty$-manifold without a boundary. (If $n=2$, then
$\mathcal K_1=\{g_1,\ g_2\}$, where $g_1,\ g_2$ are the ends of
the curves $\bar\Upsilon_1$, $\bar\Upsilon_2$.) Suppose that, in a
neighborhood of each point $g\in\mathcal K_1$, the domain $G$ is
diffeomorphic to some $n$-dimensional dihedral angle (plain angle
if $n=2$). In the domain G, we consider the nonlocal problem
\begin{gather}
 \Delta u=f_0(y)\quad (y\in G),\label{eq1Statement1}\\
 u|_{\Upsilon_i}-b_i u\big(\omega_i(y)\big)|_{\Upsilon_i}=0\quad
 (i=1,\ 2).\label{eq1Statement2}
\end{gather}
Here $b_1,\ b_2\in\mathbb R$; $\omega_i$ is an infinitely
differentiable transformation mapping some neighborhood ${\mathcal
O}_i$ of the manifold $\Upsilon_i$ onto the set $\omega({\mathcal
O}_i)$  so that $\omega_i(\Upsilon_i)\subset G$,
$\overline{\omega_i(\Upsilon_i)}\cap\partial G\ne\varnothing$, see
figures~\ref{figDomG}.a and~\ref{figDomG}.b.

Problems of type~\eqref{eq1Statement1}, \eqref{eq1Statement2} were
considered by many mathematicians (see~\cite{BitzDAN85,Kishk,GM}
and others). The most complete theory for such problems is
developed by Skubachevskii and his
pupils~\cite{SkMs86,SkDu90,SkDu91,KovSk,GurGiess,GurAsympAngle}.
In particular, Fredholm solvability of higher-order elliptic
equations with general nonlocal conditions is proved, asymptotics
for solutions near the points of conjugation of nonlocal
conditions is established, smoothness of solutions is studied. It
is shown~\cite{SkJMAA91} that the index of nonlocal problem is
equal to the index of the corresponding local one if the support
of nonlocal terms does not intersect with the points of
conjugation (see Fig.~\ref{figDomG}.a with $g_1$ and $g_2$ being
the points of conjugation). Otherwise (see Fig.~\ref{figDomG}.b),
this is not true.

Properties of nonlocal problems in bounded domains are essentially
determined by properties of model nonlocal problems in dihedral
(plain if $n=2$) angles $\Omega=\{x=(y,\ z)\in{\mathbb R}^n:
b'<\varphi<b'',\ z\in{\mathbb R}^{n-2}\}$  corresponding to the
points of conjugation of nonlocal conditions ($(\varphi,\ r)$ are
polar coordinates of~$y$). Until now~\cite{SkMs86, SkDu90,
SkDu91}, it was studied the case where the transformations
$\omega_{is}$ corresponded to linear transformations (i.e.,
compositions of rotation and expansion in $y$-plane) in model
problems. However, such a restriction is quite unnatural in
applications. Let us explain this on examples. Problem of
type~\eqref{eq1Statement1}, \eqref{eq1Statement2} is a
mathematical model for some plasma process in a bounded
domain~\cite{Sam}. Nonlocal conditions connect the plasma
temperature on the boundary of the domain with the temperature
inside the domain and at other points of the boundary.

Another important application arises in the theory of diffusion
processes. Such processes describe, for example, the Brownian
motion of a particle in the membrane $G\subset \mathbb R^n$. It is
known~\cite{Feller1,Feller2,Taira} that every diffusion process
generates some Feller semigroup. By virtue of the Hille--Iosida
theorem, the investigation of this semigroup may be reduced to the
study of an elliptic operator with boundary-value conditions
containing an integral over $\bar G$ with respect to a
non-negative Borel measure~\cite{Ventcel}. In the most difficult
case where the measure is atomic, nonlocal conditions assume the
form~\eqref{eq1Statement2}. Their probabilistic sense is as
follows: once the particle gets to a point $y\in\Upsilon_i$, it
either jumps to the point $\omega_i(y)$ with probability $b_i$
($0\le b_i\le1$) or ``dies'' with probability $1-b_i$ (in this
case, the process terminates). In general, both in the plasma
theory and theory of diffusion processes, nonlinear argument
transformations appear.

Let us mention one more application of nonlocal problems. In the
monograph~\cite{SkBook}, it is shown that in some cases a
boundary-value problem for elliptic differential-difference
equation (in particular, arising in modern aircraft technology and
modelling sandwich shells and plates~\cite{OnSk,SkBook}) can be
reduced to an elliptic equation with nonlocal conditions on shifts
of the boundary. Thus, we again obtain nonlinear transformations.
(These transformations are linear only if the boundary of domain
coincides, on certain sets, with $(n-1)$-dimensional hyperplanes.)

Other applications and references to papers devoted to nonlocal
problems can be found in~\cite{SkBook}.

\smallskip

In this paper, we consider an elliptic $2m$-order equation in a
domain $G\subset{\mathbb R}^n$ with nonlocal conditions connecting
the values of the unknown function and its derivatives on
$(n-1)$-dimensional manifolds $\Upsilon_i$ (where
$\bigcup\limits_i\bar\Upsilon_i=\partial G$) with the values on
 $\omega_{is}(\Upsilon_i)\subset G$. As we mentioned before, the
essential difficulties arise in the case where the support of
nonlocal terms
$\bigcup\limits_{i,\,s}\overline{\omega_{is}(\Upsilon_i)}$
intersects with the boundary of domain. In this situation, the
generalized solutions may have power singularities near some
set~\cite{SkMs86}. (For example, in case of
problem~\eqref{eq1Statement1}, \eqref{eq1Statement2}, these
singularities may appear near the points $g_1$ and $g_2$.)
Therefore, it is natural to consider such problems in weighted
spaces. This allows one to investigate higher-order elliptic
equations with general nonlocal conditions. We study the case
where the transformations $\omega_{is}$ correspond to nonlinear
transformations in model problems. It turns out that the problem
with nonlinear transformation is neither a small nor compact
perturbation of the corresponding local problem. Nevertheless, we
show that, when passing from linear transformations to nonlinear
ones, the operator of the problem remains Fredholm and its index
does not change.

\smallskip

Notice that a more general structure of the conjugation points and
nonlocal terms for second-order elliptic equations with nonlocal
perturbations of the Dirichlet problem was considered
in~\cite{GM}. This also justifies the importance of nonlinear
transformations~$\omega_{is}$. From our point of view, the
advantage of the approach suggested is that it allows us to study
$2m$-order elliptic equations with general boundary-value
conditions, nonlocal perturbations of which may be arbitrary
large. On the other hand, this approach also allows us to
investigate the asymptotic behavior of solutions near the
conjugation points~\cite{SkMs86,GurAsympAngle}.

\smallskip

Our paper is organized as follows. In~\S~\ref{sectStatement}, we
consider the statement of the problem and discuss the conditions
imposed on the argument transformations in nonlocal terms. Ibidem,
we introduce basic functional spaces (Sobolev spaces with a
weight) and obtain model problems in dihedral and plain angles.
In~\S~\ref{sectEx}, we give an example of nonlocal problem with
nonlinear argument transformation and show that the operator
corresponding to this problem is neither a small nor compact
perturbation of the operator corresponding to the problem with
linearized transformations. In~\S~\ref{sectOmegaK_1}, we study
some properties of nonlinear transformations near the points of
conjugation of nonlocal conditions and prove a number of lemmas
which are used in~\S~\ref{sectAprEstim} for getting a priori
estimates of solutions. In~\S~\ref{sectRightReg}, we construct a
right regularizer, which, being combined with the a priori
estimate, guarantees the Fredholm solvability of the nonlocal
problem. Finally, in~\S~\ref{sectStabInd}, we show that the index
of the problem with nonlinear argument transformations is equal to
that of the problem with the transformations linearized near the
points of conjugation of nonlocal conditions.

%% file: sect1.tex
\section{Statement of the problem in a bounded domain}\label{sectStatement}
{\bf 1.} Let $G\subset{\mathbb R}^n$ $(n\ge 2)$ be a bounded
domain with the boundary $\partial
G=\bigcup\limits_{i=1}^{N_0}\bar \Upsilon_i,$ where $\Upsilon_i$
are smooth open (in the topology of $\partial G$)
$(n-1)$-dimensional $C^\infty$-manifolds. We assume that, in a
neighborhood of each point $g\in\partial
G\setminus\bigcup\limits_{i=1}^{N_0} \Upsilon_i$, the domain $G$
is diffeomorphic to some $n$-dimensional dihedral (plain if $n=2$)
angle $\Omega=\{x=(y,\ z)\in{\mathbb R}^n: 0<b'<\varphi<b''<2\pi,\
z\in{\mathbb R}^{n-2}\}$, where $(\varphi,\ r)$ are polar
coordinates of~$y$.

We denote by ${\bf P}(x,\ D)$ and  $B_{i\mu s}(x,\ D)$
differential operators of order $2m$ and $m_{i\mu}$ respectively
with complex-valued $C^\infty$-coefficients
 ($i=1,\ \dots,\ N_0;$ $\mu=1,\
\dots,\ m;$ $s=0,\ \dots,\ S_i$). Let the operators ${\bf P}(x,\
D)$ and $B_{i\mu 0}(x,\ D)$ satisfy the following conditions (see,
for example,~\cite[Chapter~2, \S~1]{LM}).
\begin{condition}\label{condEllipPinG}
For all $x\in\bar G$, the operator ${\bf P}(x,\ D)$ is properly
elliptic.
\end{condition}
\begin{condition}\label{condComplBinG}
For all $i=1,\ \dots,\ N_0$ and $x\in\bar\Upsilon_i$, the system
$\{B_{i\mu 0}(x,\ D)\}_{\mu=1}^m$ covers the operator ${\bf P}(x,\
D)$.
\end{condition}

Let $\omega_{is}$ ($i=1,\ \dots,\ N_0;$ $s=1,\ \dots,\ S_i$) be an
infinitely differentiable transformation mapping some neighborhood
${\cal O}_i$ of the manifold $\Upsilon_i$ onto the manifold
$\omega_{is}({\cal O}_i)$ so that $\omega_{is}(\Upsilon_i)\subset
G$. We assume that the set
$$
  {\cal K}=
  \left\{\bigcup_i (\bar \Upsilon_i\setminus \Upsilon_i)\right\}\cup
  \left\{\bigcup_{i,\,s} \omega_{is}(\bar \Upsilon_i\setminus \Upsilon_i)\right\}\cup
  \left\{\bigcup_{j,\,p}\,\bigcup_{i,\,s} \omega_{jp}(\omega_{is}(\bar \Upsilon_i\setminus \Upsilon_i)\cap\Upsilon_j)\right\}
$$
can be represented in the form ${\cal K}=\bigcup\limits_{j=1}^3
{\cal K}_j,$ where
\begin{equation}\label{eqDecompK}
{\cal K}_1=\bigcup\limits_{p=1}^{N_1} {\cal K}_{1p}=
\partial G\setminus\bigcup\limits_{i=1}^{N_0}
\Upsilon_i,\quad {\cal K}_2=\bigcup\limits_{p=1}^{N_2} {\cal
K}_{2p}\subset\bigcup\limits _{i=1}^{N_0} \Upsilon_i,\quad {\cal
K}_3=\bigcup\limits_{p=1}^{N_3} {\cal K}_{3p}\subset G.
\end{equation}
Here ${\cal K}_{jp}$ are disjoint $(n-2)$-dimensional connected
$C^\infty$-manifolds without a boundary (points if $n=2$).

We consider the nonlocal boundary-value problem
\begin{equation}\label{eqPinG}
 {\bf P}(x,\ D) u=f_0(x) \quad (x\in G),
\end{equation}
\begin{equation}\label{eqBinG}
 \begin{array}{c}
   {\bf B}_{i\mu}(x,\ D)u\equiv\sum\limits_{s=0}^{S_i}(B_{i\mu s}(x,\ D)u)(\omega_{is}(x))|_{\Upsilon_i}=g_{i\mu}(x)\\
    (x\in \Upsilon_i;\ i=1,\ \dots,\ N_0;\ \mu=1,\ \dots,\ m),
 \end{array}
\end{equation}
where $(B_{i\mu s}(x,\ D)u)(\omega_{is}(x))=B_{i\mu s}(x',\
D_{x'})u(x')|_{x'=\omega_{is}(x)},$ $\omega_{i0}(x)\equiv x.$

\begin{example}
Let us consider problem~\eqref{eq1Statement1},
\eqref{eq1Statement2} in two-dimensional case, with the
transformations $\omega_i$ corresponding to Fig.~\ref{figDomGK}.
\begin{figure}[ht]
{ \hfill\epsfbox{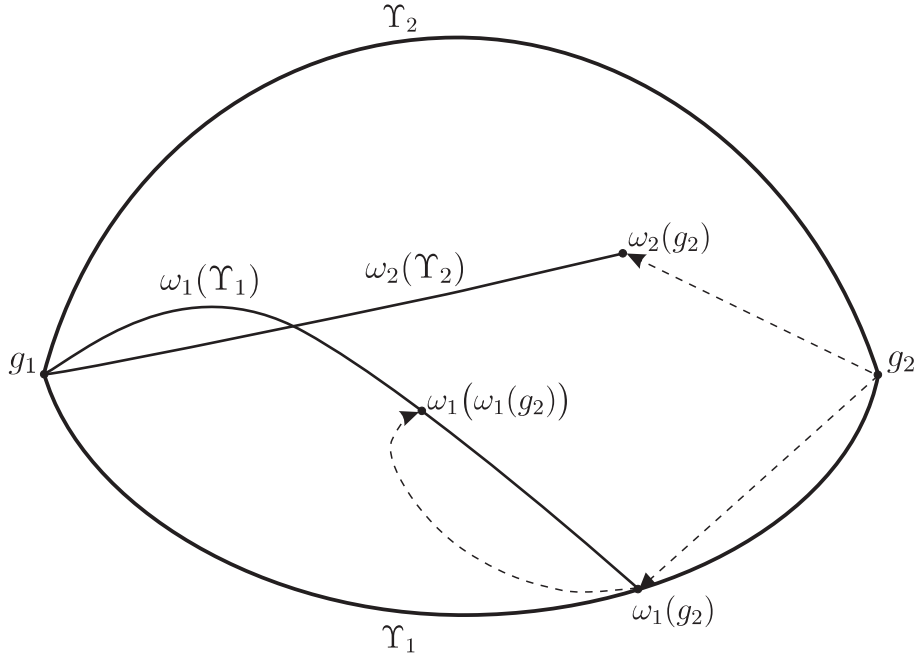}\hfill\ } \caption{The domain $G$ with
the boundary $\partial G=\bar\Upsilon_1\cup\bar\Upsilon_2$,
$n=2$.}
   \label{figDomGK}
\end{figure}
Then we have $\mathcal K_1=\{g_1,\ g_2\}$, $\mathcal
K_2=\{\omega_1(g_2)\}$, $\mathcal K_3=\big\{\omega_2(g_2),\
\omega_1\big(\omega_1(g_2)\big)\big\}$.
\end{example}

In~\cite{SkMs86}, it is shown that the solutions for
problem~(\ref{eqPinG}), (\ref{eqBinG}) may have power
singularities near the points of the set $\mathcal K_1$.
Therefore, it is natural to consider problem~(\ref{eqPinG}),
(\ref{eqBinG}) in weighted spaces. We introduce the space
$H_b^l(Q)$ as a completions of the set $C_0^\infty(\bar Q\setminus
M)$ with respect to the norm
$$
\|u\|_{H_b^l(Q)}=\left(\sum\limits_{|\alpha|\le l}\int\limits_Q
\rho^{2(b-l+|\alpha|)}|D^\alpha u|^2 dx\right)^{1/2}.
$$
Here $Q$ is the domain $G$, angle $\Omega$, or ${\mathbb R}^n$;
$M={\cal K}_1$ if $Q=G$ and $M=\{x=(y,\ z)\in{\mathbb R}^n: y=0,\
z\in{\mathbb R}^{n-2}\}$ if $Q=\Omega$ or $Q={\mathbb R}^n$;
$C_0^\infty(\bar Q\setminus M)$ is the set of infinitely
differentiable functions with compact supports being subsets of
$\bar Q\setminus M$; $l\ge 0$ is an integer; $b\in {\mathbb R};$
$\rho=\rho(x)\in C^\infty({\mathbb R}^n\setminus{\cal K}_1)$ is a
function\footnote{The existence of the function $\rho(x)$ follows
from Theorem~2~\cite[Chapter~6, \S~2]{Stein}.} satisfying $
  c_1{\rm dist}(x,\ {\cal K}_1)\le\rho(x)\le c_2{\rm dist}(x,\ {\cal K}_1)
$ ($x\in G$, $c_1,\ c_2>0$, ${\rm dist}(x,\ {\cal K}_1)$ denotes
the distance from $x$ to ${\cal K}_1$) if $Q=G$ and $\rho(x)=|y|$
if $Q=\Omega$ or $Q={\mathbb R}^n$. For $l\ge1$, we denote by
$H_b^{l-1/2}(\Upsilon)$ the space of traces on a smooth
$(n-1)$-dimensional manifold $\Upsilon\subset\bar Q$ with the norm
$$
\|\psi\|_{H_b^{l-1/2}(\Upsilon)}=\inf\|u\|_{H_b^l(Q)}\quad (u\in
H_b^l(Q):\ u|_\Upsilon=\psi).
$$

We assume that $l+2m-m_{i\mu}-1\ge0$ for all $i,\ \mu$ and
introduce the following bounded operator corresponding to nonlocal
problem~(\ref{eqPinG}), (\ref{eqBinG}):
$$
 {\bf L}=\{{\bf P}(x,\ D),\ {\bf B}_{i\mu}(x,\ D)\}: H_b^{l+2m}(G)\to
 {\cal H}_b^l(G,\ \Upsilon)=
H_b^l(G)\times\prod\limits_{i=1}^{N_0}\prod\limits_{\mu=1}^m
H_b^{l+2m-m_{i\mu}-1/2}(\Upsilon_i).
$$
From now on (unless the contrary is specified), we suppose that
$b>l+2m-1$.

Let us explain the restriction on the exponent $b$. Suppose that
the transformation $\omega_{is}$ takes a point
$g\in\bar\Upsilon_i\cap \mathcal K_1$ to the point
$\omega_{is}(g)$ so that $\omega_{is}(g)\in\mathcal K_2$ or
$\omega_{is}(g)\in\mathcal K_3$. Since the function $u(x)$ belongs
to the Sobolev space $W_2^{l+2m}$ near the point $\omega_{is}(g)$,
the function $u(\omega_{is}(x))$ belongs to the Sobolev space
$W_2^{l+2m}$ near the point $g$. However, if $b\le l+2m-1$, the
function $u(\omega_{is}(x))$ does not belong (in general) to the
weighted space $H_b^{l+2m}$. Therefore, the trace $(B_{i\mu s}(x,\
D)u)(\omega_{is}(x))|_{\Upsilon_i}$ may not belong to the weighted
space $H_b^{l+2m-m_{i\mu}-1/2}(\Upsilon_i)$, so the operator ${\bf
L}$ is not well defined. But if $b>l+2m-1$, then, by virtue of
Lemma~5.2~\cite{KovSk}, $W_2^{l+2m}(G)\subset H_b^{l+2m}(G)$.
Thus, in this case, the operator ${\bf L}$ is well defined.

Notice that, in two-dimensional case, problem~\eqref{eqPinG},
\eqref{eqBinG} can be considered in weighted spaces with arbitrary
exponent $b$ (see~\cite{SkMs86}). To this end, one should impose
some consistency conditions (generated by the transformations
$\omega_{is}$); namely, one must assume that the solutions $u$ as
well as the right-hand side $\{f_0,\ g_{i\mu}\}$ belong to the
corresponding weighted spaces not only near the set $\mathcal K_1$
but also near $\mathcal K_2$ and $\mathcal K_3$. One the one hand,
this situation is in detail considered in~\cite{SkMs86} (where the
problems with transformations linear near $\mathcal K_1$ are
studied). On the other hand, the changes described have nothing to
do with the transformations $\omega_{is}$ near $\mathcal K_1$. So,
in two-dimensional case, we will omit the proofs of corresponding
results concerning arbitrary values of $b$ (see the end
of~\S~\ref{sectRightReg}).

\medskip

{\bf 2.} Now we consider the structure of the transformations
$\omega_{is}$ near the set $\mathcal K_1$ in more detail. We
denote by $\omega_{is}^{+1}$ the transformation $\omega_{is}:{\cal
O}_i\to\omega_{is}({\cal O}_i)$ and by
$\omega_{is}^{-1}:\omega_{is}({\cal O}_i)\to{\cal O}_i$ the
transformation being inverse to $\omega_{is}$. Consider a point
$g\in{\cal K}_1$. The set of all points
$\omega_{i_ps_p}^{\pm1}(\dots\omega_{i_1s_1}^{\pm1}(g))\in{\cal
K}_1$ ($1\le s_j\le S_{i_j},\ j=1,\ \dots,\ p$) (that is, points
which can be obtained by consecutive applying to the point $g$ the
transformations $\omega_{i_js_j}^{+1}$ or $\omega_{i_js_j}^{-1}$
taking the points from ${\cal K}_1$ to those from ${\cal K}_1$) is
called an {\it orbit} of $g\in{\cal K}_1$ and denoted by
$\Orb(g)$.

We introduce the set ${\cal S}_{i1}=\{0\le s\le S_i:
\omega_{is}(\bar\Upsilon_i)\cap{\cal K}_1\ne\varnothing\}.$
Evidently, $0\in{\cal S}_{i1}.$ Let the following conditions hold.
\begin{condition}\label{condOrb}
For each $g\in{\cal K}_1$\newline (a) the set $\Orb(g)$ consists
of finitely many points $g^j$ ($j=1,\ \dots,\ N=N(g)$);\newline
(b) for the points $g^j$, there are neighborhoods
$$
 \hat{\cal V}(g^j)\subset{\cal V}(g^j)\subset{\mathbb R}^n\setminus
 \left\{\bigcup\limits_{i,\,s}\omega_{is}(\bar\Upsilon_i)\cup{\cal K}_2\cup{\cal K}_3\right\}\quad
 (s\notin {\cal S}_{i1})
$$ such that {\rm (I)}
${\cal V}(g^j)\cap{\cal V}(g^k)=\varnothing$ ($j\ne k$) and {\rm
(II)} if $g^j\in\bar\Upsilon_i$ and $\omega_{is}(g^j)=g^k,$ then
${\cal V}(g^j)\subset{\cal O}_i$ and $\omega_{is}(\hat{\cal
V}(g^j))\subset{\cal V}(g^k).$
\end{condition}

\begin{condition}\label{condCoordTrans}
For each $g\in{\cal K}_1$ and $j=1,\ \dots,\ N(g)$, there is a
non-degenerate smooth transformation $x\mapsto x'(g,\ j)$ mapping
${\cal V}(g^j)$ ($\hat{\cal V}(g^j)$) onto a neighborhood of the
origin ${\cal V}_j(0)$ ($\hat{\cal V}_j(0)$) so that\newline (a)
the images of the sets $G\cap{\cal V}(g^j)$ ($G\cap\hat{\cal
V}(g^j)$) and $\Upsilon_i\cap{\cal V}(g^j)$
($\Upsilon_i\cap\hat{\cal V}(g^j)$) are respectively the
intersection of the dihedral angle $\Omega_j=\{x=(y,\
z)\in{\mathbb R}^n:\ 0<b'_{j}<\varphi<b''_{j}<2\pi,\ z\in{\mathbb
R}^{n-2}\}$ with ${\cal V}_j(0)$ ($\hat{\cal V}_j(0)$) and the
intersection of the side of the angle $\Omega_j$ with ${\cal
V}_j(0)$ ($\hat{\cal V}_j(0)$);\newline (b) for $x\in\hat{\cal
V}(g^j)$, the transformation $\omega_{is}(x)$ ($s\in {\cal
S}_{i1}\setminus\{0\}$) in new coordinates has the form $(y',\
z')\mapsto (\omega'_{is}(y',\ z'),\ z'),$ where $\omega'_{is}(y',\
z')={\cal G}'_{is}y'+o(|x'|)$ with ${\cal G}'_{is}$ being the
operator of rotation by an angle $\varphi'_{is}$ and expansion
$\chi'_{is}>0$ times in $y'$-plane; moreover, we assume that
$\omega_{is}'(0,\ z)\equiv 0$;\newline (c) in new coordinates, the
operator ${\cal G}'_{is}$ maps the side of the corresponding angle
$\Omega_j$ ($j=j(i)$) onto an $(n-1)$-dimensional half-plane being
strictly inside an angle $\Omega_k$ ($k=k(i,\ s)$ and $j$ can be
different).
\end{condition}

Conditions~\ref{condOrb} and~\ref{condCoordTrans} are analogous to
those in~\cite{SkMs86, SkDu91}, where the transformations linear
near ${\cal K}_1$ (and arbitrary outside a neighborhood of ${\cal
K}_1$) are studied.

Condition~\ref{condOrb}~(a) is in a sense analogous to Carleman's
condition~\cite{Carleman}, which is used in the theory of nonlocal
problems with transformations mapping the boundary of domain onto
itself.

Condition~\ref{condCoordTrans}, in particular, means that if
$g\in\omega_{is}(\bar\Upsilon_i\setminus\Upsilon_i)\cap\bar\Upsilon_j\cap{\cal
K}_1\ne\varnothing$, then the surfaces
$\omega_{is}(\bar\Upsilon_i)$ and $\bar\Upsilon_j$ have different
tangent planes at the point $g$. The requirement that
$\omega_{is}'(0,\ z)\equiv 0$ is necessary for
representation~(\ref{eqDecompK}) to be possible. If
$\omega_{is}(\bar\Upsilon_i\setminus\Upsilon_i)\subset\bar
G\setminus{\cal K}_1$, then, like in~\cite{SkMs86, SkDu91}, we
have no restrictions on a geometrical structure of
$\omega_{is}(\bar\Upsilon_i)$ near $\partial G$.

\begin{remark}
One can consider the more general case where, for $x\in\hat{\cal
V}(g^j)$, the transformation $\omega_{is}(x)$ ($s\in {\cal
S}_{i1}\setminus\{0\}$) in new coordinates has the form $(y',\
z')\mapsto (\omega'_{is}(y',\ z'),\ \omega''_{is}(y',\ z')),$
where $\omega'_{is}(y',\ z')$ is the same as before,
$\omega''_{is}(y',\ z')=z'+o(|x'|)$, $\omega''_{is}(0,\ z')\equiv
z'$ (the latter guarantees that item~(a) in
Condition~\ref{condOrb} holds). However, for simplicity, we study
the transformations described in Condition~\ref{condCoordTrans}.
\end{remark}

\medskip

{\bf 3.} Let us write model problems corresponding to the points
of $\mathcal K_1$.

We fix a point $g\in{\cal K}_1.$ Let $\supp
u\subset\left(\bigcup\limits_{j=1}^{N(g)}\hat{\cal
V}(g^j)\right)\cap\bar G.$ We denote the function $u(x)$ for
$x\in{\cal V}(g^j)\cap G$ by $u_j(x)$. If $g^j\in\bar\Upsilon_i,$
$x\in\hat{\cal V}(g^j),$ $\omega_{is}(x)\in{\cal V}(g^k),$ then we
denote $u(\omega_{is}(x))$ by $u_k(\omega_{is}(x))$. Clearly,
$u(\omega_{i0}(x))\equiv u(x)\equiv u_j(x).$ Now nonlocal
problem~(\ref{eqPinG}), (\ref{eqBinG}) assumes the form
$$
 {\bf P}(x,\ D) u_j=f_0(x) \quad (x\in\hat{\cal V}(g^j)\cap G),
$$
$$
 \begin{array}{c}
   \sum\limits_{s\in{\cal S}_{i1}}(B_{i\mu s}(x,\ D)u_k)(\omega_{is}(x))|_{\Upsilon_i}=g_{i\mu}(x) \\
   (x\in \hat{\cal V}(g^j)\cap\Upsilon_i;\ i\in\{1\le i\le N_0: \hat{\cal V}(g^j)\cap\Upsilon_i\ne\varnothing\};\\
  j=1,\ \dots,\ N=N(g);\ \mu=1,\ \dots,\ m).
 \end{array}
$$

By virtue of Condition~\ref{condCoordTrans}, in new coordinates
the linear part ${\cal G}'_{is}$ of the transformation
$\omega'_{is}$ maps one of the sides of $\Omega_j$ ($j=j(i)$) onto
an $(n-1)$-dimensional half-plane being strictly inside
$\Omega_k$ ($k=k(i,\ s)$ and $j$ can be different). We denote all
these $(n-1)$-dimensional half-planes by $\Gamma_{k2},\ \dots,\
\Gamma_{k,R_k}\subset\Omega_k$. (If none of the sides of the
angles $\Omega_1,\ \dots,\ \Omega_N$ is mapped inside $\Omega_k$,
we put $R_k=1$.) We also denote $b_{k1}=b'_k$,
$b_{k,R_k+1}=b''_k$. Then the sets
$$
  \Gamma_{k\sigma}=\{x=(y,\ z)\in{\mathbb R}^n: \varphi=b_{k\sigma},\ z\in{\mathbb
  R}^{n-2}\}\quad
 (\sigma=1,\ R_k+1)
$$
are the sides of $\Omega_k$, while the half-planes $\Gamma_{kq}$
have the forms
$$
\Gamma_{kq}=\{x=(y,\ z)\in{\mathbb R}^n: \varphi=b_{kq},\
z\in{\mathbb R}^{n-2}\}\quad (q=2,\ \dots,\ R_k),
$$
where $0<b_{k1}<\cdots<b_{k,R_k+1}<2\pi$.

Let us introduce the function $U_j(x')=u_j(x(x'))$ and denote $x'$
again by $x.$ Then, by virtue of Conditions~\ref{condOrb}
and~\ref{condCoordTrans}, problem~(\ref{eqPinG}), (\ref{eqBinG})
eventually assumes the form
\begin{gather}
  {\cal P}_j(x,\ D_y,\ D_z)U_j=f_j(x) \quad
  (x\in\Omega_j),\label{eqPinOmega}\\
  \begin{array}{c}
  {\cal B}_{j\sigma\mu}(x,\ D_y,\ D_z)U\equiv B_{j\sigma\mu}(x,\ D_y,\ D_z)U_j|_{\Gamma_{j\sigma}}+\\
  +\sum\limits_{k,\,q,\,s}
                   (B_{j\sigma\mu kqs}(x,\ D_y,\ D_z)U_k)(\omega'_{j\sigma kqs}(y,\ z),\ z)|_{\Gamma_{j\sigma}}
    =g_{j\sigma\mu}(x) \quad (x\in\Gamma_{j\sigma}).
  \end{array}\label{eqBinOmega}
\end{gather}
Here (and further, until the contrary is indicated) $j,\ k=1,\
\dots,\ N;$ $\sigma=1,\ R_j+1;$ $q=2,\ \dots,\ R_k$; $\mu=1,\
\dots,\ m;$ $s=1,\ \dots,\ S_{j\sigma kq}$; ${\cal P}_j(x,\ D_y,\
D_z),$ $B_{j\sigma\mu}(x,\ D_y,\ D_z)$, and $B_{j\sigma\mu
kqs}(x,\ D_y,\ D_z)$ are operators of order $2m,$
$m_{j\sigma\mu}$, and $m_{j\sigma\mu}$ respectively with variable
$C^\infty$-coefficients; $\omega'_{j\sigma kqs}(y,\ z)={\cal
G}_{j\sigma kqs}y+o(|x|)$ with ${\cal G}_{j\sigma kqs}$ being the
operator of rotation by an angle~$\varphi_{j\sigma kq}$ and
expansion~$\chi_{j\sigma kqs}>0$ times in $y$-plane; furthermore,
$\omega'_{j\sigma kqs}(0,\ z)\equiv 0$,
$b_{k1}<b_{j\sigma}+\varphi_{j\sigma kq}=b_{kq}<b_{k,R_k+1}.$

Let us define the spaces of vector-functions:
\begin{gather*}
 H_b^{l+2m,\,N}(\Omega)=\prod_{j} H_b^{l+2m}(\Omega_j),\quad
 {\cal H}_b^{l,\,N}(\Omega,\ \Gamma)=\prod_{j} {\cal H}_b^l(\Omega_j,\
 \Gamma_j),\\
 {\cal H}_b^l(\Omega_j,\ \Gamma_j)=
 H_b^l(\Omega_j)\times\prod_{\sigma,\,\mu} H_b^{l+2m-m_{j\sigma\mu}-1/2}
   (\Gamma_{j\sigma}).
\end{gather*}
We introduce the bounded operators
$$
 {\cal L}^\omega=\{{\cal P}_j(D_y,\ D_z),\ {\cal B}^\omega_{j\sigma\mu}(D_y,\ D_z)\}:
 H_b^{l+2m,\,N}(\Omega)\to  {\cal H}_b^{l,\,N}(\Omega,\ \Gamma),
$$
$$
 {\cal L}^{\cal G}=\{{\cal P}_j(D_y,\ D_z),\ {\cal B}^{\cal G}_{j\sigma\mu}(D_y,\ D_z)\}:
 H_b^{l+2m,\,N}(\Omega)\to  {\cal H}_b^{l,\,N}(\Omega,\ \Gamma).
$$
Here\footnote{In what follows, we consider functions $U_k$ with
compact supports concentrated in a neighborhood of the origin and
such that $(\omega'_{j\sigma kqs}(y,\ z),\ z)\in\Omega_k$ for
$x\in\supp U_k$. This guarantees that the operators ${\cal
B}^\omega_{j\sigma\mu}(D_y,\ D_z)$ are well defined.}
$$
 {\cal B}^\omega_{j\sigma\mu}(D_y,\ D_z)U=B_{j\sigma\mu}(D_y,\ D_z)U_j|_{\Gamma_{j\sigma}}+
  \sum\limits_{k,\,q,\,s}
                   (B_{j\sigma\mu kqs}(D_y,\ D_z)U_k)(\omega'_{j\sigma kqs}(y,\ z),\ z)|_{\Gamma_{j\sigma}},
$$
$${\cal B}^{\cal G}_{j\sigma\mu}(D_y,\ D_z)U=B_{j\sigma\mu}(D_y,\ D_z)U_j|_{\Gamma_{j\sigma}}+
  \sum\limits_{k,\,q,\,s}(B_{j\sigma\mu kqs}(D_y,\ D_z)U_k)({\cal G}_{j\sigma kqs}y,\ z)|_{\Gamma_{j\sigma}}
$$
with ${\cal P}_j(D_y,\ D_z),$ $B_{j\sigma\mu}(D_y,\ D_z)$, and
$B_{j\sigma\mu kqs}(D_y,\ D_z)$ being the principal homogeneous
parts of the operators ${\cal P}_j(0,\ D_y,\ D_z),$
$B_{j\sigma\mu}(0,\ D_y,\ D_z)$, and $B_{j\sigma\mu kqs}(0,\ D_y,\
D_z)$ respectively.

In what follows, we will write, for short, ${\cal P}_j,$
$B_{j\sigma\mu}$, $B_{j\sigma\mu kqs}$, ${\cal
B}^\omega_{j\sigma\mu}$, and ${\cal B}^{\cal G}_{j\sigma\mu}$
instead of ${\cal P}_j(D_y,\ D_z),$ $B_{j\sigma\mu}(D_y,\ D_z)$,
$B_{j\sigma\mu kqs}(D_y,\ D_z)$, ${\cal
B}^\omega_{j\sigma\mu}(D_y,\ D_z)$, and ${\cal B}^{\cal
G}_{j\sigma\mu}(D_y,\ D_z)$ respectively.

Notice that the operator ${\cal B}^\omega_{j\sigma\mu}$ contains
nonlocal terms with nonlinear transformations $\omega'_{j\sigma
kqs}$ while the operator ${\cal B}^{\cal G}_{j\sigma\mu}$ with
linear ones ${\cal G}_{j\sigma kqs}$. Thus, the operators ${\cal
L}^\omega$ and ${\cal L}^{\cal G}$ correspond to model problems
with nonlinear and linearized transformations respectively.

As we mentioned before, the problem with transformations linear
near ${\cal K}_1$ was studied in~\cite{SkMs86,SkDu90, SkDu91}. In
particular, its Fredholm solvability was investigates.
In~\S~\ref{sectEx} of the present paper, we will show that the
operator ${\cal L}^\omega$ is neither a small nor compact
perturbation of ${\cal L}^{\cal G}$ even if the functions $U$ with
arbitrary small supports are considered. That is why, to prove the
Fredholm solvability of problem~(\ref{eqPinG}), (\ref{eqBinG})
with nonlinear transformations, we have to obtain anew a priori
estimates and construct a right regularizer
(see~\S\S~\ref{sectAprEstim}, \ref{sectRightReg}).

\medskip

{\bf 4.} Obtaining a priori estimates and constructing the right
regularizer for problem~(\ref{eqPinG}), (\ref{eqBinG}) will be
based on the invertibility of the model operators $\mathcal
L^{\mathcal G}$. Let us formulate the conditions under which the
operator $\mathcal L^{\mathcal G}$ is an isomorphism. If $n\ge3$,
then, parallel to the operator in dihedral angles, we consider a
model operator with parameter $\theta$ in plain angles. For any
angle $K=\{y\in{\mathbb R}^2:\ 0<b'<\varphi<b''<2\pi\}$, we
introduce the space $E_b^l(K)$ as a completion of $C_0^\infty(\bar
K\setminus \{0\})$ with respect to the norm
$$
 \|u\|_{E_b^l(K)}=\left(
    \sum_{|\alpha|\le l}\int\limits_K |y|^{2b}(|y|^{2(|\alpha|-l)}+1) |D_y^\alpha u(y)|^2 dy
                                       \right)^{1/2}.
$$
For $l\ge1$, we denote by $E_b^{l-1/2}(\gamma)$ the space of
traces on a ray $\gamma\subset\bar K$ with the norm
$$
 \|\psi\|_{E_b^{l-1/2}(\gamma)}=\inf\|u\|_{E_b^l(K)} \quad (u\in E_b^l(K):\  u|_\gamma = \psi).
$$
One can find the constructive definitions of the trace spaces
$H_b^{l-1/2}(\Upsilon)$ and $E_b^{l-1/2}(\gamma)$, equivalent to
the above, in~\cite[\S~1]{MP}.

We introduce the spaces of vector-functions
\begin{gather*}
 E_b^{l+2m,\,N}(K)=\prod_{j} E_b^{l+2m}(K_j),\quad
 {\cal E}_b^{l,\,N}(K,\ \gamma)=\prod_{j} {\cal E}_b^l(K_j,\
 \gamma_j),\\
 {\cal E}_b^l(K_j,\ \gamma_j)=
 E_b^l(K_j)\times\prod_{\sigma,\,\mu} E_b^{l+2m-m_{j\sigma\mu}-1/2}(\gamma_{j\sigma}),
\end{gather*}
where $K_j=\{y\in{\mathbb R}^2:\ b_{j1}<\varphi<b_{j,R_j+1}\},$
$\gamma_{j\sigma}=\{y\in{\mathbb R}^2:\ \varphi=b_{j\sigma}\}.$

We consider the bounded operator
$$
 {\cal L}^{\cal G}(\theta)=\{{\cal P}_j(D_y,\ \theta),\ {\cal B}^{\cal G}_{j\sigma\mu}(D_y,\ \theta)\}:
 E_b^{l+2m,\,N}(K)\to  {\cal E}_b^{l,\,N}(K,\ \gamma),
$$
where $\theta$ is an arbitrary point of the unit sphere
$S^{n-3}=\{\theta\in{\mathbb R}^{n-2}:\ |\theta|=1\}.$

\medskip

{\bf 5.} Let us write the operators ${\cal P}_j(D_y,\ 0),$
$B_{j\sigma\mu}(D_y,\ 0),$ $B_{j\sigma\mu kqs}(D_y,\ 0)$ in polar
coordinates: ${\cal P}_j(D_y,\ 0)=r^{-2m}\tilde{\cal
P}_j(\varphi,\ D_\varphi,\ rD_r),$ $B_{j\sigma\mu}(D_y,\
0)=r^{-m_{j\sigma\mu}}\tilde B_{j\sigma\mu}(\varphi,\ D_\varphi,\
rD_r),$ $B_{j\sigma\mu kqs}(D_y,\ 0)=r^{-m_{j\sigma\mu}}\tilde
B_{j\sigma\mu kqs}(\varphi,\ D_\varphi,\ rD_r),$ where
$D_\varphi=-i\frac{\displaystyle\partial}{\displaystyle\partial\varphi},\
D_r=-i\frac{\displaystyle\partial}{\displaystyle\partial r}.$ We
consider the analytic operator-valued function $\tilde{\cal
L}(\lambda):W_2^{l+2m,\, N}(b_1,\ b_2)\to {\cal W}_2^{l,\,
N}[b_1,\ b_2]$ given by
 \begin{multline}\notag
 \tilde{\cal L}^{\cal G}(\lambda)\tilde U=\{\tilde{\cal P}_j(\varphi,\ D_\varphi,\ \lambda)\tilde U_j,\
  \tilde B_{j\sigma\mu}(\varphi,\ D_\varphi,\ \lambda)\tilde U_j(\varphi)|_{\varphi=b_{j\sigma}}+\\
           + \sum\limits_{k,\,q,\,s}  e^{(i\lambda-m_{j\sigma\mu})\ln\chi_{j\sigma kqs}}{\tilde B}_{j\sigma\mu kqs}(\varphi,\ D_\varphi,\ \lambda)
              \tilde U_k(\varphi+\varphi_{j\sigma kqs})|_{\varphi=b_{j\sigma}}\},
 \end{multline}
where
\begin{gather*}
 W_2^{l+2m,\, N}(b_1,\ b_2)=\prod_{j} W_2^{l+2m}(b_{j1},\
 b_{j,R_j+1}),\quad
 {\cal W}_2^{l,\, N}[b_1,\ b_2]=\prod_{j} {\cal W}_2^{l}[b_{j1},\
 b_{j,R_j+1}],\\
 {\cal W}_2^{l}[b_{j1},\ b_{j,R_j+1}]=W_2^l(b_{j1},\ b_{j,R_j+1}) \times{\mathbb C}^{2m}.
\end{gather*}

By virtue of Lemmas~2.1, 2.2~\cite{SkDu90}, there exists a
finite-meromorphic operator-valued function $(\tilde{\cal L}^{\cal
G})^{-1}(\lambda)$ such that $(\tilde{\cal L}^{\cal
G})^{-1}(\lambda)$ is the inverse to $\tilde{\cal L}^{\cal
G}(\lambda)$ if $\lambda$ is not a pole of $(\tilde{\cal L}^{\cal
G})^{-1}(\lambda)$; furthermore, for every pole $\lambda_0$, there
is a $\delta>0$ such that the set $\{\lambda\in{\mathbb
C}:0<|\Im\lambda-\Im\lambda_0|<\delta\}$ contains no poles of
$(\tilde{\cal L}^{\cal G})^{-1}(\lambda)$.

If $n=2$, then, by Theorem~2.1~\cite{SkDu90}, the operator ${\cal
L}^{\cal G}$ is an isomorphism if and only if the line
$\Im\lambda=b+1-l-2m$ contains no poles of $(\tilde{\cal L}^{\cal
G})^{-1}(\lambda)$.

Suppose that $n\ge3$ and assume that the system
$\{B_{j\sigma\mu}(D_y,\ D_z)\}_{\mu=1}^m$ is normal on
$\Gamma_{j\sigma}$ and the orders $m_{j\sigma\mu}$ of the
operators $B_{j\sigma\mu}(D_y,\ D_z)$, $B_{j\sigma\mu kqs}(D_y,\
D_z)$ are less or equal to $2m-1$. In this case, by virtue of
Theorem~9.1 \cite{GurGiess}, the operator ${\cal L}^{\cal
G}(\theta)$ is Fredholm if and only if the line
$\Im\lambda=b+1-l-2m$ contains no poles of $(\tilde{\cal L}^{\cal
G})^{-1}(\lambda)$. By Theorem~3.3~\cite{SkDu90}, if, in addition,
$\dim\ker( {\cal L}^{\cal G}(\theta))=\codim{\cal R}({\cal
L}^{\cal G}(\theta))=0$ for $b$ replaced by $b-l$, $l$ replaced by
$0$, and all $\theta\in S^{n-3},$ then the operator
$$
 {\cal L}^{\cal G}=\{{\cal P}_j(D_y,\ D_z),\ {\cal B}^{\cal G}_{j\sigma\mu}(D_y,\ D_z)\}:
 H_b^{l+2m,\,N}(\Omega)\to  {\cal H}_b^{l,\,N}(\Omega,\ \Gamma)
$$
is an isomorphism (see the corresponding example
in~\cite[\S~10]{GurGiess}). Notice that if ${\cal L}^{\cal G}$ is
not an isomorphism, then ${\cal L}^{\cal G}(\theta)$ is not
Fredholm (see Theorem~9.3~\cite{GurGiess}).

Since the operators ${\cal L}^\omega,$ ${\cal L}^{\cal G},$ ${\cal
L}^{\cal G}(\theta)$, and $\tilde{\cal L}^{\cal G}(\lambda)$
corresponding to problem~(\ref{eqPinOmega}), (\ref{eqBinOmega})
depend on the choice of $g\in {\cal K}_1$, we denote them by
${\cal L}^\omega_g,$ ${\cal L}^{\cal G}_g,$ ${\cal L}^{\cal
G}_g(\theta)$, and $\tilde{\cal L}^{\cal G}_g(\lambda)$
respectively.

%% file: sect2.tex
\section{Example of nonlocal problem with nonlinear
argument transformations}\label{sectEx} In this section, we show
on a simple example that a problem with a transformation nonlinear
in a neighborhood of ${\cal K}_1$ is neither a small nor compact
perturbation of the problem with the linearized transformation.

\begin{figure}[ht]
{ \hfill\epsfbox{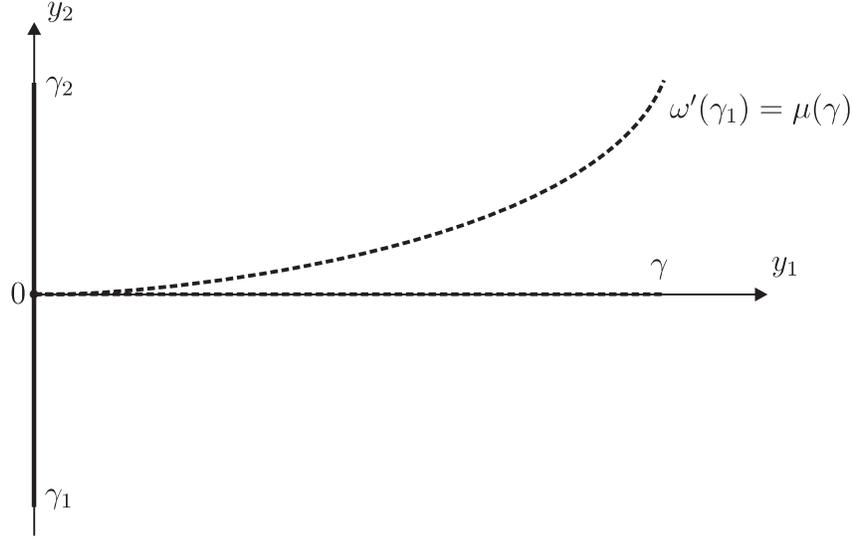}\hfill\ } \caption{The angle $K$ of
opening $\pi$.}
   \label{figAngleK}
\end{figure}

{\bf 1.} Let us assume for simplicity that problem~(\ref{eqPinG}),
(\ref{eqBinG}) is considered in a plain domain. Let the model
problem~(\ref{eqPinOmega}), (\ref{eqBinOmega}) corresponding to
some point of ${\cal K}_1$ have the form
\begin{gather*}
 \triangle u=f(y)\qquad (y\in K),\\
\begin{aligned}
 u|_{\gamma_1}+u(\omega'(y))|_{\gamma_1}&=g_1(y)\quad (y\in\gamma_1),\\
 u|_{\gamma_2}&=g_2(y)\quad (y\in\gamma_2).
 \end{aligned}
\end{gather*}
Here $K=\{y\in{\mathbb R}^2: r>0,\ |\varphi|<\pi/2\}$ is a plain
angle (of opening $\pi$) with the sides $\gamma_i=\{y\in{\mathbb
R}^2: r>0,\ \varphi=(-1)^i \pi/2\}$ ($i=1,\ 2$). We suppose that
$\omega'(y)=\mu({\cal G}y)$, where ${\cal G}$ is the operator of
rotation by the angle $\pi/2$ mapping $\gamma_1$ onto a ray
$\gamma=\{y\in{\mathbb R}^2: r>0,\ \varphi=0\}$;
$$
 \mu: (y_1,\ y_2)\mapsto \left(\frac{y_1}{\sqrt{1+y_1^2}},\ y_2+\frac{y_1^2}{\sqrt{1+y_1^2}}\right)
$$
is an infinitely differentiable transformation mapping $\gamma$
onto the curve $\mu(\gamma)$, which is tangent to $\gamma$ at the
origin (see Fig.~\ref{figAngleK}).

The operators ${\cal L}^\omega,\ {\cal L}^{\cal G}:
H_b^{l+2}(K)\to H_b^l(K)\times\prod\limits_{i=1}^2
H_b^{l+3/2}(\gamma_i)$ corresponding to the model problems with
nonlinear and linearized transformations have the form
$$
{\cal L}^\omega u=\{\triangle u,\ u|_{\gamma_1}+u(\omega'(y))|_{\gamma_1},\
u|_{\gamma_2}\},
$$
$$
{\cal L}^{\cal G} u=\{\triangle u,\ u|_{\gamma_1}+u({\cal G} y)|_{\gamma_1},\
u|_{\gamma_2}\}.
$$

Clearly, a non-zero component of the difference ${\cal L}^{\cal
G}u-{\cal L}^\omega u$ is
$$
  u({\cal G}y)|_{\gamma_1} - u(\omega'(y))|_{\gamma_1}=u(y)|_\gamma - u(\mu(y))|_\gamma.
$$

We introduce the operator $A_\varepsilon: H_b^{l+2}(K)\to
H_b^{l+3/2}(\gamma)$ with the domain $D(A_\varepsilon)=\{u\in
H_b^{l+2}(K) : \supp u\subset \{r<\varepsilon\}\cap\bar K\}$ given
by
$$
 A_\varepsilon u(y)=u(y)|_\gamma - u(\mu(y))|_\gamma.
$$

In this example, we prove that one cannot make the operator
$A_\varepsilon$ small or compact, choosing sufficiently small
$\varepsilon$. For simplicity, we show this in the case where
$A_\varepsilon$ acts from $H_b^{1}(K)$ to $H_b^{1/2}(\gamma)$. The
general case can be considered in the same way. We shall construct
a sequence $u_\varepsilon\in D(A_\varepsilon)$, $\varepsilon\to
0$, such that
$$
 \|u_\varepsilon|_\gamma - u_\varepsilon(\mu (\cdot))|_\gamma\|_{H_b^{1/2}(\gamma)}\ge c\|u_\varepsilon\|_{H_b^1(K)},
$$
where $c>0$ is independent of $\varepsilon$.

Let us write the restriction of $\mu$ on $\gamma$ in polar
coordinates $(\varphi,\ r)$:
$$
 \mu|_\gamma: (0,\ r)\mapsto (\Phi(r),\ r),
$$
where $\Phi(r)=\arctg r$. Clearly, $\Phi(0)=0$, $\Phi(1)=\pi/4$,
$\frac{1}{\sqrt{2}}\le\frac{\Phi}{r},\frac{d\Phi}{dr}\le1$ on
$[0,1]$.

Let us consider the transformation
$$
 \tilde\mu: (\varphi,\ r)\mapsto (\varphi+\Phi(r),\ r).
$$
One can see that $u(\mu(y))|_\gamma=u(\tilde\mu(y))|_\gamma$ since
$\mu|_\gamma=\tilde\mu|_\gamma$. Therefore, without loss of
generality, we may assume that the transformation $\mu$ is given
by
$$
 \mu: (\varphi,\ r)\mapsto (\varphi+\Phi(r),\ r).
$$

Notice that the norm of any function $u\in H_b^1(K)$ written in
polar coordinates is equivalent to
$$
 \left(\sum\limits_{|\alpha|\le1}\int\limits_{0}^\infty\int\limits_{-\pi/2}^{\pi/2}r^{2b-1}
          |(rD_r)^{\alpha_1}D_\varphi^{\alpha_2}u(\varphi,\ r)|^2\,d\varphi dr\right)^{1/2}.
$$

Set $r=e^{-t}$; then, in new coordinates, the transformation $\mu$
assumes the form
$$
 \mu: (\varphi,\ t)\mapsto (\varphi+\Phi(e^{-t}),\ t).
$$
Putting $v(\varphi,\ t)=u(\varphi,\ e^{-t})$, we see that the norm
$\|u\|_{H_b^1(K)}$ is equivalent to the norm
\begin{equation}\label{eqNormQ}
 \|v\|_{W_{2,b}^1(Q)}= \left(\sum\limits_{|\alpha|\le1}\int\limits_{-\infty}^\infty\int\limits_{-\pi/2}^{\pi/2}e^{-2bt}
          |D_t^{\alpha_1}D_\varphi^{\alpha_2}v(\varphi,\ t)|^2\,d\varphi dt\right)^{1/2},
\end{equation}
where $Q=\{t\in{\mathbb R},\ |\varphi|<\pi/2\}$ and $W_{2,b}^1(Q)$
is the space with norm~(\ref{eqNormQ}). Evidently, $W_{2,0}^1(Q)$
coincides with the Sobolev space $W_{2}^1(Q)$.

Since the norms $\|v\|_{W_{2,b}^1(Q)}$ and
$\|e^{-bt}v\|_{W_{2}^1(Q)}$ are equivalent, it suffices to study
the case where $b=0$. In what follows, we consider functions
$v(\varphi,\ t)$ with the support being a subset of the strip
$\{|\varphi|<\pi/2\}$. Putting $v=0$ for $|\varphi|\ge\pi/2$, we
obtain $\|v\|_{W_2^1(Q)}=\|v\|_{W_{2}^1({\mathbb R}^2)}$.

Thus, our task is reduced to constructing a sequence $v_s\in
W_2^1({\mathbb R}^2)$ such that $\supp v_s\subset\{t>2s,\
|\varphi|<\pi/2\}$ and
$$
 \|v_s(0,\ t) - v_s(\Phi(e^{-t}),\ t)\|_{W_2^{1/2}({\mathbb R})}\ge c\|v_s\|_{W_{2}^1({\mathbb R}^2)},
$$
where $c>0$ is independent of $s$.

To this end, we pass from variables $(\varphi,\ t)$ to $(\varphi,\
\tau)$: we introduce the sets
$$
 Q_s=\left\{|\theta|\le\frac{\pi}{2},\  2s\le\tau\le 2s+1\right\},\quad s=0,\ 1,\ 2,\
 \dots,
$$
and put
\begin{equation}\label{eqThetaTau}
\varphi=F(\theta,\ \tau),\quad t=\tau.
\end{equation}
Here $F(\theta,\ \tau)=\theta e^{2s}\Phi(e^{-\tau})$ for
$(\theta,\ \tau)\in Q_s$, $s=0,\ 1,\ 2,\ \dots$, and $F(\theta,\
\tau)$ is extended onto ${\mathbb R}^2\setminus
\bigcup\limits_{s=0}^\infty Q_s$ so that the
transformation~(\ref{eqThetaTau}) remains continuously
differentiable with the Jacobian $\frac{\partial F}{\partial
\theta}$ such that
\begin{equation}\label{eqJacobian}
 0<c_1\le \left|\frac{\partial F}{\partial \theta}\right|\le c_2\quad \text{on } {\mathbb
R}^2.
\end{equation}
Such an extension does exist: indeed,
$$
 \frac{\partial F}{\partial \theta}=e^{2s}\Phi(e^{-\tau}),\quad
\frac{\partial F}{\partial \tau}=-\theta e^{-\tau+2s}
\frac{d\Phi}{d r}\Big|_{r=e^{-\tau}},\qquad (\theta,\ \tau)\in
Q_s;
$$
therefore (by virtue of the above properties of $\Phi$), in
$\bigcup\limits_{s=0}^\infty Q_s$ the function $F(\theta,\ \tau)$
is continuously differentiable with respect to $\theta$ and $\tau$
and inequalities~\eqref{eqJacobian} hold.

One easily sees that, under change of
variables~(\ref{eqThetaTau}), the segment $Q_s\cap\{\theta=0\}$ is
an image of the corresponding segment of the line $\{\varphi=0\}$.
Furthermore, the transformation $\mu$ on $Q_s$ has the form
\begin{equation}\label{eqOmegaInQs}
 \mu: (\theta,\ \tau)\mapsto (\theta+e^{-2s},\ \tau),\quad (\theta,\ \tau)\in Q_s.
\end{equation}

\begin{figure}[ht]
{ \hfill\epsfbox{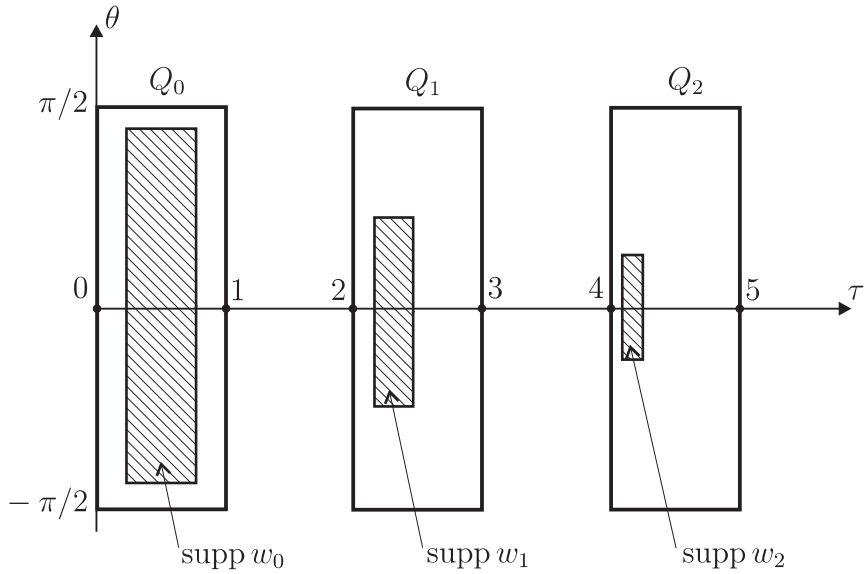}\hfill\ }
    \caption{The supports of $w_s$ are contained in the hatched domains.}
   \label{figSuppW_s}
\end{figure}

We consider functions $f,\ g\in C^\infty({\mathbb R})$ such that
$\supp f\subset\{|\theta|<\frac{\pi}{2}\}$, $f(0)\ne f(1)$, $\supp
g\subset \{0<\tau<1\}$, $g(\tau)\not\equiv0$ and define the
sequence $w_s(\theta,\ \tau)=f_s(\theta)g_s(\tau)$, where
$$
 f_s(\theta)=f(\theta e^{2s}),\quad g_s(\tau)=g((\tau-2s)e^{2s}),\quad s=0,\ 1,\ 2,\ \dots
$$
Clearly, $\supp w_s\subset Q_s$ (see Fig.~\ref{figSuppW_s}).

We have
\begin{multline}\label{eqNormR2}
 \|w_s\|^2_{W_2^1({\mathbb R}^2)}=\|f_s\|^2_{L_2({\mathbb R})}\|g_s\|^2_{L_2({\mathbb R})}+
\left\|\frac{d f_s}{d\theta}\right\|^2_{L_2({\mathbb
R})}\|g_s\|^2_{L_2({\mathbb R})}+
\|f_s\|^2_{L_2({\mathbb R})}\left\|\frac{d g_s}{d\tau}\right\|^2_{L_2({\mathbb R})}=\\
=e^{-4s}\|f\|^2_{L_2({\mathbb R})}\|g\|^2_{L_2({\mathbb R})}+
\left\|\frac{d f}{d\theta}\right\|^2_{L_2({\mathbb
R})}\|g\|^2_{L_2({\mathbb R})}+ \|f\|^2_{L_2({\mathbb
R})}\left\|\frac{d g}{d\tau}\right\|^2_{L_2({\mathbb R})}.
\end{multline}
Analogously, using the fact that the norm in $W_2^{1/2}({\mathbb
R})$ is given by
$$
 \|g\|_{W_2^{1/2}({\mathbb R})}=\left(\|g\|^2_{L_2({\mathbb R})}+
 \int\limits_{-\infty}^\infty\int\limits_{-\infty}^\infty\frac{|g(\tau_1)-g(\tau_2)|^2}{|\tau_1-\tau_2|^2}d\tau_1 d\tau_2\right)^{1/2}
$$
(see~\cite{Slob}) and the form~(\ref{eqOmegaInQs}) of the
transformation $\mu$ in coordinates $(\theta,\ \tau)$, we get
\begin{multline}\label{eqNormR1}
  \|w_s|_{\theta=0}-w_s(\mu (\cdot))|_{\theta=0}\|^2_{W_2^{1/2}({\mathbb R})}=
|f_s(0)-f_s(e^{-2s})|^2 \|g_s\|^2_{W_2^{1/2}({\mathbb R})}\ge\\
\ge
|f(0)-f(1)|^2\int\limits_{-\infty}^\infty\int\limits_{-\infty}^\infty\frac{|g(\tau_1)-g(\tau_2)|^2}{|\tau_1-\tau_2|^2}d\tau_1
d\tau_2.
\end{multline}
From~(\ref{eqNormR2}) and (\ref{eqNormR1}), it follows that
$$
 \|w_s|_{\theta=0}-w_s(\mu (\cdot))|_{\theta=0}\|^2_{W_2^{1/2}({\mathbb R})}\ge c\|w_s\|^2_{W_2^1({\mathbb R}^2)}.
$$

\medskip

{\bf 2.} Using the sequence $w_s$, one can easily show that, for
any $\varepsilon$, the operator $A_\varepsilon$ is not compact.
Indeed, the sequence $w_s$ is bounded in $W_2^1({\mathbb R}^2)$.
However, one cannot choose from $w_s|_{\theta=0}-w_s(\mu
(\cdot))|_{\theta=0}$ a subsequence convergent in
$W_2^{1/2}({\mathbb R})$, since, according to~(\ref{eqNormR1}),
for all natural $s\ne h$ the expression
 \begin{multline}\notag
  \|[w_s|_{\theta=0}-w_s(\mu (\cdot))|_{\theta=0}]-[w_h|_{\theta=0}-w_h(\mu(\cdot))|_{\theta=0}]\|_{W_2^{1/2}({\mathbb R})}=\\
  =\|w_s|_{\theta=0}-w_s(\mu(\cdot))|_{\theta=0}\|_{W_2^{1/2}({\mathbb R})}+
  \|w_h|_{\theta=0}-w_h(\mu(\cdot))|_{\theta=0}\|_{W_2^{1/2}({\mathbb R})}
\end{multline}
is bounded from below by a positive constant.

%% file: sect3.tex
\section{Argument transformations near the set ${\cal K}_1$}\label{sectOmegaK_1}
From the results of~\S~\ref{sectEx}, it follows that, to prove the
Fredholm solvability of the problem with transformations nonlinear
near ${\cal K}_1$, one has to obtain anew  a priori estimates and
construct the right regularizer. To this end, we start by studying
some properties of the transformations $\omega_{is}$ near the set
${\cal K}_1$.

We fix a point $g\in{\cal K}_1$, make, for each $j=1,\ \dots,\
N=N(g)$, the change of variables $x\mapsto x'(g,\ j)$, and
consider the transformations $\omega'_{j\sigma kqs}(y,\ z)$ for
$(y,\ z)\in{\cal V}_{\varepsilon_0}(0)=\{x\in {\mathbb R}^n:
|x|<\varepsilon_0\}.$ The number $\varepsilon_0$ is supposed to be
small so that $\overline{{\cal
V}_{\varepsilon_0}(0)}\subset\hat{\cal V}_j(0),$ $j=1,\ \dots,\
N$. In the sequel, we shall impose some additional conditions on
$\varepsilon_0$.

{\bf 1.} Before we proceed to study the transformations
$\omega_{is}$, let us prove an auxiliary result, which will be
used for proving a lemma on a representation of $\omega_{is}$ in
polar coordinates (see Lemma~\ref{lOmega}).
\begin{lemma}\label{lh/r}
Let $h=h(r,\ z)$ be a function such that $|D^k_rD_z^\alpha h|\le
c_{k\alpha}$ for $r\ge0$, $z\in{\mathbb R}^{n-2}$,
 $(r^2+|z|^2)^{1/2}\le\varepsilon_0$. Set $f(r,\ z)=r^{-l}h(r,\ z)$
for some $l\in{\mathbb N}$ and assume that $|f|\le c$. Then
$|D_r^kf|\le c_k$ for $r\ge0$, $z\in{\mathbb R}^{n-2}$,
$(r^2+|z|^2)^{1/2}\le\varepsilon_0$, and any $k=1,\ 2\dots$
\end{lemma}
\begin{proof}
1) First, we consider the case where $l=1$, that is $f(r,\
z)=r^{-1}h(r,\ z)$. By Leibnitz' formula, we have
$$
 \frac{\displaystyle \partial^k f(r,\ z)}{\displaystyle\partial r^k}=\sum\limits_{s=0}^k\frac{\displaystyle(-1)^s k!}{\displaystyle (k-s)!}
 r^{-s-1}\frac{\displaystyle \partial^{k-s} h(r,\ z)}{\displaystyle\partial r^{k-s}}.
$$
Expanding $\frac{\displaystyle \partial^{k-s}
h}{\displaystyle\partial r^{k-s}}$ by the Taylor formula near
$r=0$ and using the boundedness of the derivatives of $h$, we
obtain
 \begin{multline}\notag
  \frac{\displaystyle \partial^k f(r,\ z)}{\displaystyle\partial r^k}=\sum\limits_{s=0}^k\frac{\displaystyle(-1)^s k!}{\displaystyle (k-s)!}
  r^{-s-1}\left[\sum\limits_{p=0}^s\frac{\displaystyle 1}{\displaystyle p!}
  \frac{\displaystyle \partial^{k-s+p} h}{\displaystyle\partial r^{k-s+p}}(0,\ z)r^p+
  \frac{\displaystyle \partial^{k+1} h}{\displaystyle\partial r^{k+1}}(\varkappa_{rz} r,\ z)r^{s+1}\right]=\\
  =\sum\limits_{s=0}^k\sum\limits_{p=0}^s \frac{\displaystyle(-1)^s k!}{\displaystyle (k-s)!p!}
  \frac{\displaystyle\partial^{k-s+p} h}{\displaystyle\partial r^{k-s+p}}(0,\
  z)r^{-s-1+p}+O(1),
 \end{multline}
where $\varkappa_{rz}\in(0,\ 1)$.

Putting $p'=s-p$ in the last sum and denoting $p'$ again by $p$,
we get
$$
 \frac{\displaystyle \partial^k f(r,\ z)}{\displaystyle\partial r^k}=\sum\limits_{s=0}^k\sum\limits_{p=0}^s
 \frac{\displaystyle(-1)^s k!}{\displaystyle (k-s)!(s-p)!}
 \frac{\displaystyle\partial^{k-p} h}{\displaystyle\partial r^{k-p}}(0,\ z)r^{-p-1}+O(1).
$$
Write the coefficient $a_p(z)$ at $r^{-p-1}$ on the right-hand
side of the last identity:
 \begin{multline}\notag
 a_p(z)=\frac{\displaystyle\partial^{k-p} h}{\displaystyle\partial r^{k-p}}(0,\ z)
 \sum\limits_{s=p}^k\frac{\displaystyle(-1)^s k!}{\displaystyle (k-s)!(s-p)!}=\\
 =\frac{\displaystyle\partial^{k-p} h}{\displaystyle\partial r^{k-p}}(0,\ z)(-1)^p
  \sum\limits_{s=0}^{k-p}k(k-1)\cdot\dots\cdot(k-(s+p)+1)\frac{\displaystyle 1}{\displaystyle s!}(-1)^s,\quad p=0,\ \dots,\ k.
 \end{multline}
Since $|r^{-1}h(r,\ z)|\le c$ by assumption, we have $h(0,\
z)\equiv 0$; therefore, $a_k(z)\equiv 0$. On the other hand,
notice that, for $0\le p<k$, we have
 \begin{multline}\notag
 0= \frac{\displaystyle d^p}{\displaystyle dt^p}(t+1)^k\Big|_{t=-1}=
 \left.\left(\sum\limits_{s=0}^{k-p}k(k-1)\cdot\dots\cdot(k-(s+p)+1)\frac{\displaystyle 1}{\displaystyle s!}t^s\right)\right|_{t=-1}=\\
 =\sum\limits_{s=0}^{k-p}k(k-1)\cdot\dots\cdot(k-(s+p)+1)\frac{\displaystyle 1}{\displaystyle s!}(-1)^s.
 \end{multline}
Thus, $a_p(z)\equiv 0$ for all $p=0,\ \dots,\ k$, and the lemma is
proved for $l=1$.

2) For $l\ge2$, we use the mathematical induction method. Let the
lemma be true for $l=1,\ \dots, l_1-1$. We claim that it is true
for $l=l_1$. We have $f=r^{-1}f_1$, where $f_1=r^{-(l_1-1)}h$.
Since $|f|\le c$, it follows that $|f_1|\le c$, and, therefore, by
the inductive assumption (for $l=l_1-1$) the estimate
$|D^k_rD_z^\alpha f_1|\le c_{k\alpha}$ holds. Applying the
inductive assumption once more (now, for $l=1$), we get the
conclusion of the lemma for $r^{-1}f_1$, that is, for
$f=r^{-l_1}h$.
\end{proof}

Now let us proceed to investigate the transformations
$\omega_{is}$. The following lemma describes the structure of
$\omega'_{j\sigma kqs}$ in cylindrical coordinates. Such a
representation turns out to be convenient for the study of
nonlocal problems in weighted spaces.

\begin{lemma}\label{lOmega}
For sufficiently small $\varepsilon_0$, the transformation
$\omega'_{j\sigma kqs}(y,\ z)|_{\Gamma_{j\sigma}\cap{\cal
V}_{\varepsilon_0}(0)}$ can be represented in polar coordinates in
the form
 \begin{equation}\label{eqOmega01}
  (b_{j\sigma},\ r)\mapsto\big(b_{kq}+\Phi_{j\sigma kqs}(r,\ z),\ \chi_{j\sigma kqs}r+R_{j\sigma kqs}(r,\ z)\big)\quad \mbox{for }
   (r^2+|z|^2)^{1/2}\le\varepsilon_0,
 \end{equation}
where $\Phi_{j\sigma kqs}(r,\ z),\ R_{j\sigma kqs}(r,\ z)$ are
infinitely differentiable functions such that
 \begin{equation}\label{eqOmega02}
  |\Phi_{j\sigma kqs}|\le c \varepsilon_0,\quad |R_{j\sigma kqs}|\le c \varepsilon_0 r,
 \end{equation}
 \begin{equation}\label{eqOmega03}
  |D_r^k D_z^\alpha\Phi_{j\sigma kqs}|\le c_{k\alpha},\quad
  |D_r^k D_z^\alpha(R_{j\sigma kqs}/r)|\le c_{k\alpha}.
 \end{equation}
Here $k+|\alpha|\ge1;$
 $c,\ c_{k\alpha}>0$ are independent of $\varepsilon_0.$
\end{lemma}
\begin{proof}
Let $\omega'_{j\sigma kqs}(y,\ z)=(\omega^1_{j\sigma kqs}(y,\ z),\
\omega^2_{j\sigma kqs}(y,\ z)).$ By
condition~\ref{condCoordTrans}, we have $\omega^i_{j\sigma
kqs}(0,\ z)\equiv 0$ ($i=1,\ 2$); therefore, the Teylor formula
near $r=0$ implies
\begin{equation}\label{eqOmega1}
  \omega^i_{j\sigma kqs}(r\cos b_{j\sigma},\ r\sin b_{j\sigma},\ z)=\left(\frac{\partial \omega^i_{j\sigma kqs}}
      {\partial y_1}(0,\ z)\cos b_{j\sigma}+\frac{\partial \omega^i_{j\sigma kqs}}
      {\partial y_2}(0,\ z)\sin b_{j\sigma}\right)r+O(r^2).
\end{equation}
Here $O(r^2)$ is a function with absolute values majorized by
$cr^2$, where $c$ is independent of $r$ and $z$. (To verify this,
one should write the remainder of the Teylor formula in Lagrange's
form and use smoothness of $\omega^i_{j\sigma kqs}$.) Expanding
$\frac{\partial \omega^i_{j\sigma kqs}}{\partial y_1}(0,\ z)$ and
$\frac{\partial \omega^i_{j\sigma kqs}}{\partial y_2}(0,\ z)$ by
the Teylor formula near $z=0$, from~(\ref{eqOmega1}) we obtain
\begin{equation}\label{eqOmega2}
      \omega^i_{j\sigma kqs}=\left(\frac{\partial \omega^i_{j\sigma kqs}}
      {\partial y_1}(0)\cos b_{j\sigma}+\frac{\partial \omega^i_{j\sigma kqs}}
      {\partial y_2}(0)\sin b_{j\sigma}\right)r+O(|z|) r+O(r^2).
\end{equation}

Notice that $\frac{\partial \omega^1_{j\sigma kqs}}{\partial
y_1}(0)\cos b_{j\sigma}+\frac{\partial \omega^1_{j\sigma kqs}}
{\partial y_2}(0)\sin b_{j\sigma}$ and $\frac{\partial
\omega^2_{j\sigma kqs}}{\partial y_1}(0)\cos
b_{j\sigma}+\frac{\partial \omega^2_{j\sigma kqs}} {\partial
y_2}(0)\sin b_{j\sigma}$ are not simultaneously equal to zero.
(This follows from non-degeneracy of the Jacobian of the
transformation $(y,\ z)\mapsto (\omega'_{j\sigma kqs}(y,\ z),\ z)$
at the origin.) For definiteness, we assume that
\begin{equation}\label{eqOmega3}
 \frac{\partial \omega^1_{j\sigma kqs}}{\partial y_1}(0)\cos b_{j\sigma}+\frac{\partial \omega^1_{j\sigma kqs}}
 {\partial y_2}(0)\sin b_{j\sigma}\ne 0.
\end{equation}
Hence, by virtue of~(\ref{eqOmega2}),
\begin{equation}\label{eqOmega4}
 \omega^1_{j\sigma kqs}\ne0\quad \mbox{for } (r^2+|z|^2)^{1/2}\le\varepsilon_0
\end{equation}
with $\varepsilon_0$ small enough, and the transformation
$\omega'_{j\sigma kqs}|_{\Gamma_{j\sigma}\cap{\cal
V}_{\varepsilon_0}(0)}$ in polar coordinates has the form
\begin{equation}\label{eqOmega5}
 \left(b_{j\sigma},\ r\right)\mapsto
 \left(\arctg\frac{\omega^2_{j\sigma kqs}}
     {\omega^1_{j\sigma kqs}}+\pi l,\
     \sqrt{\sum\limits_{i=1}^2(\omega^i_{j\sigma kqs})^2}\right),
\end{equation}
where $l=0$ if $\omega^1_{j\sigma kqs}>0$ and $\omega^2_{j\sigma
kqs}\ge0$, $l=1$ if $\omega^1_{j\sigma kqs}<0$, $l=2$ if
$\omega^1_{j\sigma kqs}>0$ and $\omega^2_{j\sigma kqs}<0.$

From~(\ref{eqOmega2}) and the Teylor formula, it follows that
$$
 \begin{array}{c}
 \arctg\frac{\omega^2_{j\sigma kqs}}
     {\omega^1_{j\sigma kqs}}=
 \arctg\frac{\frac{\partial \omega^2_{j\sigma kqs}}{\partial y_1}(0)\cos b_{j\sigma}+\frac{\partial \omega^2_{j\sigma kqs}}
{\partial y_2}(0)\sin b_{j\sigma}}{\frac{\partial \omega^1_{j\sigma kqs}}
{\partial y_1}(0)\cos b_{j\sigma}+\frac{\partial \omega^1_{j\sigma kqs}}{\partial y_2}(0)\sin b_{j\sigma}}+O(|z|)+O(r),\\
\\
\sqrt{\sum\limits_{i=1}^2 (\omega^i_{j\sigma kqs})^2}
 =r\sqrt{ \sum\limits_{i=1}^2\left(\frac{\partial \omega^i_{j\sigma kqs}}{\partial y_1}(0)\cos b_{j\sigma}
     +\frac{\partial \omega^i_{j\sigma kqs}}{\partial y_2}(0)\sin b_{j\sigma}\right)^2}
 +O(|z|)r+O(r^2).
\end{array}
$$
Setting
$$
  b_{kq}=\arctg\frac{\frac{\partial \omega^2_{j\sigma kqs}}{\partial y_1}(0)\cos b_{j\sigma}+\frac{\partial \omega^2_{j\sigma kqs}}
{\partial y_2}(0)\sin b_{j\sigma}}{\frac{\partial \omega^1_{j\sigma kqs}}
{\partial y_1}(0)\cos b_{j\sigma}+\frac{\partial \omega^1_{j\sigma
kqs}}{\partial y_2}(0)\sin b_{j\sigma}}+\pi l,
$$
$$
 \chi_{j\sigma kqs}=\sqrt{ \sum\limits_{i=1}^2\left(\frac{\partial \omega^i_{j\sigma kqs}}{\partial y_1}(0)\cos b_{j\sigma}+
      \frac{\partial \omega^i_{j\sigma kqs}}{\partial y_2}(0)\sin b_{j\sigma}\right)^2},
$$
we get formula~(\ref{eqOmega01}) and
inequalities~(\ref{eqOmega02}).

Let us prove the first inequality in~(\ref{eqOmega03}).
By~(\ref{eqOmega4}), we have $\left|\frac{\omega^2_{j\sigma kqs}}
     {\omega^1_{j\sigma kqs}}\right|\le c$ for $(r^2+|z|^2)^{1/2}\le\varepsilon_0$.
Therefore, by virtue of~(\ref{eqOmega01}) and~(\ref{eqOmega5}), it
suffices to prove that the derivatives
$D_r^kD_z^\alpha\frac{\omega^2_{j\sigma kqs}} {\omega^1_{j\sigma
kqs}}$ are bounded. Clearly, we have
$$
 \frac{\omega^2_{j\sigma kqs}}
     {\omega^1_{j\sigma kqs}}=
\frac{r^{-1}\omega^2_{j\sigma kqs}}
     {r^{-1}\omega^1_{j\sigma kqs}}.
$$
From~(\ref{eqOmega2}) and~(\ref{eqOmega3}), it follows that
$r^{-1}\omega^1_{j\sigma kqs}\ne0$ for
$(r^2+|z|^2)^{1/2}\le\varepsilon_0$. Hence, it suffices to prove
that
$$
 |D_r^kD_z^\alpha(r^{-1}\omega^i_{j\sigma kqs})|
 =|D_r^k(r^{-1}D_z^\alpha \omega^i_{j\sigma kqs})|\le c_{k\alpha},\quad i=1,\ 2.
$$
But the function $D_z^\alpha \omega^i_{j\sigma kqs}$ is infinitely
differentiable for $(r^2+|z|^2)^{1/2}\le\varepsilon_0$;
furthermore, since $\omega^i_{j\sigma kqs}(0,\ z)\equiv 0$, we
have $D_z^\alpha \omega^i_{j\sigma kqs}=O(r)$. Therefore, $
|r^{-1}D_z^\alpha \omega^i_{j\sigma kqs}|\le c_\alpha$. Now the
conclusion of the lemma follows from Lemma~\ref{lh/r}.

Similarly, one can prove the second inequality
in~(\ref{eqOmega03}). From~(\ref{eqOmega01}) and~(\ref{eqOmega5}),
it follows that
$$
\frac{R_{j\sigma kqs}(r,\ z)}{r}=\sqrt{\sum\limits_{i=1}^2
 \frac{(\omega^i_{j\sigma kqs})^2}{r^2}}
 -\chi_{j\sigma kqs}.
$$
By virtue of~(\ref{eqOmega2}) and~(\ref{eqOmega3}), we obtain
$\sum\limits_{i=1}^2(\omega^i_{j\sigma kqs})^2/r^2\ne 0$ for
$(r^2+|z|^2)^{1/2}\le\varepsilon_0$; therefore, it suffices to
prove that
$$
 \left|D_r^kD_z^\alpha\sum\limits_{i=1}^2(\omega^i_{j\sigma kqs})^2/r^2
  \right|\le c_{k\alpha}.
$$
But the function $D_z^\alpha\sum\limits_{i=1}^2(\omega^i_{j\sigma
kqs})^2$ is infinitely differentiable for
$(r^2+|z|^2)^{1/2}\le\varepsilon_0$; furthermore, since
$\omega^i_{j\sigma kqs}(0,\ z)\equiv 0$, we have $D_z^\alpha
\sum\limits_{i=1}^2(\omega^i_{j\sigma kqs})^2=O(r^2)$. Hence,
$\left|D_z^\alpha\sum\limits_{i=1}^2(\omega^i_{j\sigma
kqs})^2/r^2\right|\le
 c_\alpha$, and the conclusion of the lemma again follows from Lemma~\ref{lh/r}.
\end{proof}

\medskip

{\bf 2.} Denote $\delta=\min\{b_{j,q+1}-b_{jq}\}/2\ (j=1,\ \dots,\
N;\ q=1,\ \dots,\ R_j),$ $d_1=\min\{1,\ \chi_{j\sigma kqs}\}/2,$
$d_2=2\max\{1,\ \chi_{j\sigma kqs}\}.$ Let $\varepsilon_0$ be so
small that
\begin{equation}\label{eqPhiLessDelta}
|\Phi_{j\sigma kqs}|\le\delta/2,\quad |R_{j\sigma
kqs}|\le\chi_{j\sigma kqs}r/2\quad \mbox{for }
(r^2+|z|^2)^{1/2}\le\varepsilon_0/d_1.
\end{equation}
The existence of such an $\varepsilon_0$ follows from
Lemma~\ref{lOmega}.

We introduce infinitely differentiable functions
$\zeta_{j\sigma,i}(\varphi),$ $\zeta_{kq,i}(\varphi)$ such that
\begin{equation}\label{eqZetaJSigma}
 \begin{array}{c}
 \zeta_{j\sigma,i}(\varphi)=1\ \mbox{for }
 |b_{j\sigma}-\varphi|\le\delta/2^{i+1},\quad
     \zeta_{j\sigma,i}(\varphi)=0\ \mbox{for } |b_{j\sigma}-\varphi|\ge\delta/2^{i},\\
     \\
\zeta_{kq,i}(\varphi)=\zeta_{j\sigma,i}(\varphi-\varphi_{j\sigma
kq}),
 \end{array}
\end{equation}
$i=0,\ \dots,\ 4$. Clearly, $\zeta_{kq,i}(\varphi)=1\ \mbox{for }
|b_{kq}-\varphi|\le\delta/2^{i+1},\
     \zeta_{kq,i}(\varphi)=0\ \mbox{for } |b_{kq}-\varphi|\ge\delta/2^{i}.$

Let us consider the transformation $\tilde\omega'_{j\sigma
kqs}(y,\ z)$ that are given in polar coordinates by
\begin{equation}\label{eqHatOmega'}
 (\varphi,\ r)\mapsto(\varphi+\varphi_{j\sigma kq}+\Phi_{j\sigma kqs}(r,\ z),\ \chi_{j\sigma kqs}r+R_{j\sigma kqs}(r,\ z)).
\end{equation}
By virtue of Lemma~\ref{lOmega}, we have $\tilde\omega'_{j\sigma
kqs}(y,\ z)|_{\Gamma_{j\sigma}\cap {\cal V}_{\varepsilon_0}(0)}=
\omega'_{j\sigma kqs}(y,\ z)|_{\Gamma_{j\sigma}\cap {\cal
V}_{\varepsilon_0}(0)}$; therefore, in what follows, we can assume
that the transformation $\omega'_{j\sigma kqs}(y,\ z)$ is given
by~(\ref{eqHatOmega'}). Notice that now $\omega'_{j\sigma kqs}(y,\
z)$ may have (in general) a singularity at the origin, since the
new transformation $\omega'_{j\sigma kqs}(y,\ z)$ coincides with
the old one $\omega'_{j\sigma kqs}(y,\ z)$ only on
$\Gamma_{j\sigma}\cap {\cal V}_{\varepsilon_0}(0)$.

For any function $W(y,\ z)$, we denote $\hat W(y,\
z)=W(\omega'_{j\sigma kqs}({\cal G}_{j\sigma kqs}^{-1}y,\ z),\
z)$. By virtue of Lemma~\ref{lOmega}, $\omega'_{j\sigma kqs}({\cal
G}_{j\sigma kqs}^{-1}y,\ z)$ in polar coordinates has the form
\begin{equation}\label{eqHatOmegaG'}
 (\varphi,\ r)\mapsto(\varphi+\Phi'_{j\sigma kqs}(r,\ z),\ r+R'_{j\sigma kqs}(r,\ z)),
\end{equation}
where $\Phi'_{j\sigma kqs}(r,\ z)=\Phi_{j\sigma kqs}(\chi_{j\sigma
kqs}^{-1}r,\ z)$, $R'_{j\sigma kqs}(r,\ z)=R_{j\sigma
kqs}(\chi_{j\sigma kqs}^{-1}r,\ z)$. It is easy to see that
$\Phi'_{j\sigma kqs}$ and $R'_{j\sigma kqs}$ also satisfy
inequalities~(\ref{eqOmega02}), (\ref{eqOmega03}).

\begin{lemma}\label{lOmega'H}
For sufficiently small $\varepsilon_0$ and any $W\in
H_b^l(\Omega_k)$ with $\supp W\subset\bar \Omega_k\cap{\cal
V}_{\varepsilon_0}(0)$ we have $\zeta_{kq,1}\hat W\in
H_b^l(\Omega_k)$ and
  $$
   \|\zeta_{kq,1}\hat W\|_{H_b^l(\Omega_k)}\le c\|W\|_{H_b^l(\Omega_k)},
  $$
where $q=2,\ \dots,\ R_k$; $c>0$ is independent of $W$ and
$\varepsilon_0$.
\end{lemma}
\begin{proof}
In the proof, we shall use the following obvious assertion:
\begin{equation}\label{eqWinHEquiv}
 W\in H_b^l(\Omega_k)\Longleftrightarrow D^\alpha W\in H_{b+|\alpha|-l}^0(\Omega_k),\ |\alpha|\le l.
\end{equation}
From formula~(\ref{eqHatOmegaG'}) and
inequalities~(\ref{eqPhiLessDelta}), it follows that the
transformation~(\ref{eqHatOmegaG'}) maps $\overline{{\cal
V}_{\varepsilon_0}(0)}\cap\{x:|\varphi-b_{kq}|<\delta\}\cap
\Omega_k$ into $\Omega_k$ for $q=2,\ \dots,\ R_k$. Furthermore,
inequalities~(\ref{eqOmega02}) and (\ref{eqOmega03}) imply that,
for small $\varepsilon_0$, the absolute value of the Jacobian of
transformation~(\ref{eqHatOmegaG'}) is bounded and does not vanish
in $\overline{{\cal
V}_{\varepsilon_0}(0)}\cap\{x:|\varphi-b_{kq}|<\delta\}\cap\Omega_k$.
This proves the lemma for $l=0$ and $\zeta_{kq,0}$ substituted for
$\zeta_{kq,1}$.

Let us consider functions $\zeta_{kq,0}^{p}\in C_0^\infty({\mathbb
R})$ ($p=0,\ \dots,\ l$) such that
$\zeta_{kq,0}^{0}=\zeta_{kq,0}$, $\zeta_{kq,0}^{l}=\zeta_{kq,1}$,
and $\zeta_{kq,0}^{p-1}(\varphi)=1$ for
$\varphi\in\supp\zeta_{kq,0}^{p}$ ($p=1,\ \dots,\ l$). Let us
assume that the lemma is true for $l=p-1$ and $\zeta_{kq,0}^{p-1}$
substituted for $\zeta_{kq,1}$. We claim that it is true for $l=p$
and $\zeta_{kq,0}^{p}$ substituted for $\zeta_{kq,1}$ ($p\ge1$).
Indeed, let $W\in H_b^{p}(\Omega_k)$; then $\frac{\displaystyle
1}{\displaystyle r} \frac{\displaystyle \partial W}{\displaystyle
\partial  \varphi}$,
$\frac{\displaystyle \partial W}{\displaystyle  \partial r}$,
$\frac{\displaystyle \partial W}{\displaystyle  \partial z_\xi}\in
H_b^{p-1}(\Omega_k)$, $\xi=1,\ \dots,\ n-2$. Therefore, by the
inductive assumption, we have
$\zeta_{kq,0}^{p-1}\widehat{\frac{\displaystyle 1}{\displaystyle
r} \frac{\displaystyle \partial W}{\displaystyle \partial
\varphi}}$, $\zeta_{kq,0}^{p-1}\widehat{\frac{\displaystyle
\partial W}{\displaystyle  \partial r}}$,
$\zeta_{kq,0}^{p-1}\widehat{\frac{\displaystyle \partial
W}{\displaystyle  \partial z_\xi}}\in H_b^{p-1}(\Omega_k)$. From
this, relations
 \begin{equation}\label{eqComplexDif}
 \begin{array}{c}
 \frac{\displaystyle 1}{\displaystyle r}
      \frac{\displaystyle \partial \hat W_k}{\displaystyle \partial  \varphi}=
   \widehat{\frac{\displaystyle 1}{\displaystyle r} \frac{\displaystyle \partial W}{\displaystyle \partial  \varphi}}
   \cdot(1+\frac{\displaystyle R'_{j\sigma qks}}{\displaystyle r}),\\
\\
      \frac{\displaystyle \partial \hat W_k}{\displaystyle \partial  r}=
  \widehat{\frac{\displaystyle 1}{\displaystyle r} \frac{\displaystyle \partial W}{\displaystyle \partial  \varphi}}
   \cdot(1+\frac{\displaystyle R'_{j\sigma kqs}}{\displaystyle  r})\cdot
   r\frac{\displaystyle\partial \Phi'_{j\sigma kqs}}{\displaystyle\partial r}
 +\widehat{\frac{\displaystyle \partial W}{\displaystyle  \partial r}}
   \cdot(1+\frac{\displaystyle\partial R'_{j\sigma kqs}}{\displaystyle\partial r}),\\
\\
      \frac{\displaystyle \partial \hat W_k}{\displaystyle \partial  z_\xi}=
   \widehat{\frac{\displaystyle 1}{\displaystyle r} \frac{\displaystyle \partial W}{\displaystyle \partial  \varphi}}
   \cdot(1+\frac{\displaystyle R'_{j\sigma kqs}}{\displaystyle  r})\cdot
   r\frac{\displaystyle\partial \Phi'_{j\sigma kqs}}{\displaystyle\partial z_\xi}
 +\widehat{\frac{\displaystyle \partial W}{\displaystyle  \partial r}}
   \cdot\frac{\displaystyle\partial R'_{j\sigma kqs}}{\displaystyle\partial z_\xi}
 +\widehat{\frac{\displaystyle \partial W}{\displaystyle  \partial z_\xi}},
 \end{array}
\end{equation}
inequalities~(\ref{eqOmega02}), (\ref{eqOmega03}), and
Lemma\footnote { Lemma~2.1~\cite{KondrTMMO67} (and Lemmas~2.2,
3.5, 3.6~\cite{KondrTMMO67}, see below) is proved by Kondrat'ev
for domains with angular or conical points. However, it is easy to
see that it remains true for the domains with edges under
consideration.}~2.1~\cite{KondrTMMO67}, we get
\begin{equation}\label{eqWinHEquiv1}
 \zeta_{kq,0}^{p-1}\frac{\displaystyle 1}{\displaystyle r} \frac{\displaystyle \partial \hat W}{\displaystyle \partial  \varphi},\
 \zeta_{kq,0}^{p-1}\frac{\displaystyle \partial \hat W}{\displaystyle  \partial r},\
 \zeta_{kq,0}^{p-1}\frac{\displaystyle \partial \hat W}{\displaystyle  \partial z_\xi}\in H_b^{p-1}(\Omega_k).
\end{equation}
Furthermore, the relation $W\in H_b^{p}(\Omega_k)$, embedding
$H_b^{p}(\Omega_k)\subset H_{b-p}^0(\Omega_k)$, and the conclusion
of the lemma for $l=0$ imply $\zeta_{kq,0}^{p}\hat W\in
H_{b-p}^0(\Omega_k)$. From this, (\ref{eqWinHEquiv}),
and~(\ref{eqWinHEquiv1}), it follows that
$D^\alpha(\zeta_{kq,0}^{p}\hat W)\in
H_{b+|\alpha|-p}^0(\Omega_k)$, $|\alpha|\le p$. Once more
using~(\ref{eqWinHEquiv}), we complete the proof.
\end{proof}

Thus, we proved that the operator $W\mapsto\zeta_{kq,1} \hat W$ is
bounded in $H_b^l(\Omega_k)$.

\begin{lemma}\label{lGammaOmega'}
For any $W\in H_b^{l}(\Omega_k)$ with $\supp W\subset\bar
\Omega_k\cap{\cal V}_{\varepsilon_0}(0)$ and any multi-index
$\gamma,$ $1\le|\gamma|\le l$, the following inequality holds:
 \begin{equation}\label{eqlGammaOmega'}
  \|\zeta_{kq,2}D^\gamma \hat W-\zeta_{kq,2}\widehat{D^\gamma W}\|_{H_b^{l-|\gamma|}(\Omega_k)}
  \le c\varepsilon_0\|W\|_{H_b^{l}(\Omega_k)},
 \end{equation}
where $q=2,\ \dots,\ R_k$; $c>0$ is independent of $W$ and
$\varepsilon_0$.
\end{lemma}
\begin{proof}
We introduce functions $\zeta_{kq,1}^{p}\in C_0^\infty({\mathbb
R})$ ($p=1,\ \dots,\ l$) such that
$\zeta_{kq,1}^{1}=\zeta_{kq,1}$, $\zeta_{kq,1}^{l}=\zeta_{kq,2}$,
and $\zeta_{kq,1}^{p-1}(\varphi)=1$ for
$\varphi\in\supp\zeta_{kq,1}^{p}$ ($p=2,\ \dots,\ l$).

Let $|\gamma|=1$; then it suffices to prove
inequality~(\ref{eqlGammaOmega'}) for the case where the operator
$D^\gamma$ is replaced by $\frac{\displaystyle 1}{\displaystyle r}
\frac{\displaystyle
\partial}{\displaystyle \partial  \varphi},$ $\frac{\displaystyle
\partial}{\displaystyle  \partial r},$ $\frac{\displaystyle
\partial}{\displaystyle  \partial z_\xi}$. Let us consider
the operator $\frac{\displaystyle 1}{\displaystyle r}
\frac{\displaystyle
\partial}{\displaystyle \partial  \varphi}$ (the other operators can be considered
in the same way). Combining the first relation
in~(\ref{eqComplexDif}) with Leibniz' formula, we get
 \begin{multline}\notag
\left\|\zeta_{kq,1}^{1}\frac{\displaystyle 1}{\displaystyle r}
\frac{\displaystyle \partial\hat W}{\displaystyle \partial
\varphi}- \zeta_{kq,1}^{1}\widehat{\frac{\displaystyle
1}{\displaystyle r} \frac{\displaystyle \partial W}{\displaystyle
\partial  \varphi}}\right\|_
{H_b^{l-1}(\Omega_k)}^2=\left\|\zeta_{kq,1}^{1}
   \widehat{\frac{\displaystyle 1}{\displaystyle r} \frac{\displaystyle \partial W}{\displaystyle \partial  \varphi}}
     \frac{\displaystyle R'_{j\sigma qks}}{\displaystyle r}\right\|_{H_b^{l-1}(\Omega_k)}^2\le\\
 \le k_1\sum\limits_{|\alpha|\le l-1}\sum\limits_{|\beta|\le|\alpha|}
 \int\limits_{\Omega_k}r^{2(b+|\alpha|-(l-1))}
 \left|D^{\alpha-\beta}\frac{\displaystyle R'_{j\sigma qks}}{\displaystyle r}\right|^2  \left|D^{\beta}(\zeta_{kq,1}^{1}
     \widehat{\frac{\displaystyle 1}{\displaystyle r} \frac{\displaystyle \partial W}{\displaystyle \partial
     \varphi}})\right|^2\,dx.
 \end{multline}
From this, the last inequality in~(\ref{eqOmega02}), and the last
inequality in~(\ref{eqOmega03}), we obtain
\begin{equation}\label{eqGammaOmega'1}
 \left\|\zeta_{kq,1}^{1}\frac{\displaystyle 1}{\displaystyle r} \frac{\displaystyle \partial\hat W}{\displaystyle \partial  \varphi}-
 \zeta_{kq,1}^{1}\widehat{\frac{\displaystyle 1}{\displaystyle r} \frac{\displaystyle \partial W}{\displaystyle \partial  \varphi}}\right\|_
  {H_b^{l-1}(\Omega_k)}^2\le k_2\varepsilon_0^2\left\|\zeta_{kq,1}^{1}
     \widehat{\frac{\displaystyle 1}{\displaystyle r} \frac{\displaystyle \partial W}{\displaystyle \partial  \varphi}}\right\|
       _{H_b^{l-1}(\Omega_k)}^2.
\end{equation}
Estimate~(\ref{eqGammaOmega'1}) and Lemma~\ref{lOmega'H} prove the
lemma for $|\gamma|=1$ and $\zeta_{kq,1}^{1}$ substituted for
$\zeta_{kq,2}$.

We assume that the lemma is true for $1\le|\gamma|\le p-1$ and
$\zeta_{kq,1}^{p-1}$ substituted for $\zeta_{kq,2}$. Let us prove
that it is true for $|\gamma|=p$ and $\zeta_{kq,1}^{p}$
substituted for $\zeta_{kq,2}$ ($p\ge2$). We have
\begin{multline}\label{eqGammaOmega'2}
 \|\zeta_{kq,1}^{p}D^\gamma \hat W-\zeta_{kq,1}^p\widehat{D^\gamma W}\|_{H_b^{l-|\gamma|}(\Omega_k)}\le
 \|\zeta_{kq,1}^{p}D^{|\gamma|-1}(D^1\hat W)-\zeta_{kq,1}^pD^{|\gamma|-1}\widehat{D^1 W}\|_{H_b^{l-|\gamma|}(\Omega_k)}+\\
+\|\zeta_{kq,1}^{p}D^{|\gamma|-1}\widehat{D^1
W}-\zeta_{kq,1}^p\widehat{D^{|\gamma|-1}(D^1
W)}\|_{H_b^{l-|\gamma|}(\Omega_k)}\le
k_3(\|\zeta_{kq,1}^{p-1}D^1\hat W-\zeta_{kq,1}^{p-1}\widehat{D^1 W}\|_{H_b^{l-1}(\Omega_k)}+\\
+\|\zeta_{kq,1}^{p}D^{|\gamma|-1}\widehat{D^1
W}-\zeta_{kq,1}^p\widehat{D^{|\gamma|-1}(D^1
W)}\|_{H_b^{l-|\gamma|}(\Omega_k)}),
\end{multline}
where $D^{|\gamma|-1}$ and $D^1$ are some derivatives of order
$|\gamma|-1$ and $1$ respectively. By the inductive assumption,
for each of the two norms on the right-hand side
of~(\ref{eqGammaOmega'2}), the following estimates hold:
$$
 \|\zeta_{kq,1}^{p-1}D^1\hat W-\zeta_{kq,1}^{p-1}\widehat{D^1 W}\|_{H_b^{l-1}(\Omega_k)}\le k_4
 \varepsilon_0\|W\|_{H_b^{l}(\Omega_k)},
$$
$$
 \|\zeta_{kq,1}^{p}D^{|\gamma|-1}\widehat{D^1 W}-\zeta_{kq,1}^p\widehat{D^{|\gamma|-1}(D^1 W)}\|_{H_b^{l-|\gamma|}(\Omega_k)}\le
 k_5  \varepsilon_0\|D^1 W\|_{H_b^{l-1}(\Omega_k)}
\le k_6 \varepsilon_0\|W\|_{H_b^{l}(\Omega_k)}.
$$
This and~(\ref{eqGammaOmega'2}) imply the conclusion of the lemma.
\end{proof}

Notice that the multiplier $\varepsilon_0$ appears
in~(\ref{eqlGammaOmega'}) since the minuend and subtrahend both
contain the same transformation $\omega'_{j\sigma kqs}({\cal
G}_{j\sigma kqs}^{-1}y,\ z)$, but the minuend is the derivative
$D^\gamma$ of the transformed function $\hat W$ while the
subtrahend is the transformation of the derivative $D^\gamma W$.

\begin{lemma}\label{lBOmegaG}
For any $U_k\in H_b^{l+2m}(\Omega_k)$ with $\supp U_k\subset\bar
\Omega_k\cap{\cal V}_{\varepsilon_0}(0),$ the following inequality
holds:
 \begin{multline}\label{eqBOmegaG}
  \|(B_{j\sigma\mu kqs}U_k)({\cal G}_{j\sigma kqs}y,\ z)|_{\Gamma_{j\sigma}}-
      (B_{j\sigma\mu kqs}U_k)(\omega'_{j\sigma kqs}(y,\ z),\ z)|_{\Gamma_{j\sigma}}\|_
      {H_b^{l+2m-m_{j\sigma\mu}-1/2}(\Gamma_{j\sigma})}\le \\
  \le c(\varepsilon_0\|U_k\|_{H_b^{l+2m}(\Omega_k)}+
  \|\zeta_{kq,3}U_k-\zeta_{kq,3}\hat U_k\|_{H_b^{l+2m}(\Omega_k)}),
 \end{multline}
where $c>0$ is independent of $U$ and $\varepsilon_0$.
\end{lemma}
\begin{proof}
Using the boundedness of the trace operator in weighted spaces, we
get
\begin{multline}\label{eqBOmegaG1}
  \|(B_{j\sigma\mu kqs}U_k)({\cal G}_{j\sigma kqs}y,\ z)|_{\Gamma_{j\sigma}}-
    (B_{j\sigma\mu kqs}U_k)(\omega'_{j\sigma kqs}(y,\ z),\ z)|_{\Gamma_{j\sigma}}\|_
    {H_b^{l+2m-m_{j\sigma\mu}-1/2}(\Gamma_{j\sigma})}\le\\
  \le k_1\|\zeta_{kq,4}B_{j\sigma\mu kqs}U_k-\zeta_{kq,4}\widehat {B_{j\sigma\mu kqs}U_k}\|_
     {H_b^{l+2m-m_{j\sigma\mu}}(\Omega_k)}\le \\
 \le k_1(\|\zeta_{kq,4}B_{j\sigma\mu kqs}U_k-\zeta_{kq,4}B_{j\sigma\mu kqs}\hat U_k\|_
   {H_b^{l+2m-m_{j\sigma\mu}}(\Omega_k)}+\\
 +\|\zeta_{kq,4}B_{j\sigma\mu kqs}\hat U_k-\zeta_{kq,4}\widehat {B_{j\sigma\mu kqs}U_k}\|_
   {H_b^{l+2m-m_{j\sigma\mu}}(\Omega_k)}).
\end{multline}
Let us estimate the first norm on the right-hand side
of~(\ref{eqBOmegaG1}) as follows:
\begin{equation}\label{eqBOmegaG2}
 \|\zeta_{kq,4}B_{j\sigma\mu kqs}U_k-\zeta_{kq,4}B_{j\sigma\mu kqs}\hat U_k\|_
   {H_b^{l+2m-m_{j\sigma\mu}}(\Omega_k)}\le k_2\|\zeta_{kq,3}U_k-\zeta_{kq,3}\hat U_k\|_{H_b^{l+2m}(\Omega_k)}.
\end{equation}
The second norm on the right-hand side of~(\ref{eqBOmegaG1}) can
be estimated with the help of Lemma~\ref{lGammaOmega'}:
\begin{equation}\label{eqBOmegaG3}
 \|\zeta_{kq,4}B_{j\sigma\mu kqs}\hat U_k-\zeta_{kq,4}\widehat {B_{j\sigma\mu kqs}U_k}\|_
   {H_b^{l+2m-m_{j\sigma\mu}}(\Omega_k)})\le k_3\varepsilon_0\|U_k\|_{H_b^{l+2m}(\Omega_k)}.
\end{equation}
From~(\ref{eqBOmegaG1})--(\ref{eqBOmegaG3}), the conclusion of the
lemma follows.
\end{proof}

Notice that the right-hand side of~(\ref{eqBOmegaG}) contains the
norm of the difference of the non-transformed function and the
transformed one. To estimate such differences, we need the
following result.

\begin{lemma}\label{lW-W}
For any $W\in H_{b+1}^{1}(\Omega_k)$ with $\supp W\subset\bar
\Omega_k\cap{\cal V}_{\varepsilon_0}(0)$, the following inequality
holds:
 \begin{equation}\label{eqW-W}
  \|\zeta_{kq,1}W-\zeta_{kq,1}\hat W\|_{H_{b}^{0}(\Omega_k)}\le c\varepsilon_0\|W\|_{H_{b+1}^{1}(\Omega_k)}.
 \end{equation}
where $c>0$ is independent of $W$ and $\varepsilon_0$.
\end{lemma}
\begin{proof}
Writing the arguments of the functions $W$ and $\hat W$ in
cylindrical coordinates, we obtain
\begin{multline}\label{eqW-W1}
 \|\zeta_{kq,1}W-\zeta_{kq,1}\hat W\|_{H_{b}^{0}(\Omega_k)}\le
 \|\zeta_{kq,1}W(\varphi,\ r,\ z)-\zeta_{kq,1}W(\varphi+\Phi'_{j\sigma kqs}(r,\ z),\ r,\ z)\|_{H_{b}^{0}(\Omega_k)}+\\
 +\|\zeta_{kq,1}W(\varphi+\Phi'_{j\sigma kqs}(r,\ z),\ r,\ z)-
   \zeta_{kq,1}W(\varphi+\Phi'_{j\sigma kqs}(r,\ z),\ r+R'_{j\sigma kqs}(r,\ z),\ z)\|_{H_{b}^{0}(\Omega_k)}.
\end{multline}
Using the Schwartz inequality, we estimate the square of the first
norm on the right-hand side of~(\ref{eqW-W1}):
 \begin{multline}\notag
 \|\zeta_{kq,1}W(\varphi,\ r,\ z)-\zeta_{kq,1}W(\varphi+\Phi'_{j\sigma kqs}(r,\ z),\ r,\ z)\|_{H_{b}^{0}(\Omega_k)}^2=\\
=\int\limits_{{\mathbb R}^{n-2}}dz\int\limits_0^\infty
r^{2b}r\,dr\int\limits_{b_{k1}}^{b_{k2}}
\left|\zeta_{kq,1}\int\limits_{\varphi}^{\varphi+\Phi'_{j\sigma kqs}(r,\ z)}
  \frac{\displaystyle \partial W}{\displaystyle \partial  \varphi'}\,d\varphi'\right|^2\,d\varphi\le\\
\int\limits_{{\mathbb R}^{n-2}}dz\int\limits_0^\infty
r^{2b}r\,dr\int\limits_{b_{k1}}^{b_{k2}} |\zeta_{kq,1}|^2 |\Phi'_{j\sigma
kqs}(r,\ z)|\cdot \left|\int\limits_{\varphi}^{\varphi+\Phi'_{j\sigma kqs}(r,\
z)}
  \left|\frac{\displaystyle \partial W}{\displaystyle \partial  \varphi'}\right|^2\,d\varphi'\right|\,d\varphi.
 \end{multline}
Taking into account the restrictions on the support of the
functions $W$ and $\zeta_{kq,1}$ and
inequalities~(\ref{eqPhiLessDelta}), we can change the order of
integration with respect to $\varphi$ and $\varphi'$; as a result,
using~(\ref{eqOmega02}), we get
 \begin{multline}\notag
 \|\zeta_{kq,1}W(\varphi,\ r,\ z)-\zeta_{kq,1}W(\varphi+\Phi'_{j\sigma kqs}(r,\ z),\ r,\ z)\|_{H_{b}^{0}(\Omega_k)}^2\le\\
 \le k_1 \int\limits_{{\mathbb R}^{n-2}}dz\int\limits_0^\infty r^{2b}r|\Phi'_{j\sigma kqs}(r,\ z)|^2\,dr
 \int\limits_{b_{k1}}^{b_{k2}}\left|\frac{\displaystyle \partial W}{\displaystyle \partial  \varphi}\right|^2\,d\varphi\le\\
 \le k_2\varepsilon_0^2\int\limits_{{\mathbb R}^{n-2}}dz\int\limits_0^\infty r^{2(b+1)}r\,dr
 \int\limits_{b_{k1}}^{b_{k2}}\left|\frac{\displaystyle 1}{\displaystyle r}
 \frac{\displaystyle \partial W}{\displaystyle \partial  \varphi}\right|^2\,
 d\varphi\le k_3\varepsilon_0^2\|W\|_{H_{b+1}^{1}(\Omega_k)}^2.
 \end{multline}
Similarly, one can estimate the square of the second norm on the
right-hand side of~(\ref{eqW-W1}).
\end{proof}

Thus, the multiplier $\varepsilon_0$ appears in~(\ref{eqW-W}) if
one increases the order of differentiation by 1. (The left-hand
side of~(\ref{eqW-W}) contains the norm in $H_b^0(\Omega_k)$ while
the right-hand side does in $H_{b+1}^1(\Omega_k)$.) This can be
explained as follows: unlike in~(\ref{eqlGammaOmega'}), in this
case one estimates the difference of the two functions the first
one of which does not contain a transformation while the second
one does.

%% file: sect4.tex
\section{A priori estimates of solutions}\label{sectAprEstim}
In this section, we prove an a priori estimate for the operator
${\bf L}$, which guarantees that its kernel is of finite dimension
and its range is closed.

{\bf 1.} First, we prove an a priori estimate for functions with
the support being a subset of some neighborhood of ${\cal K}_1$.
To this end, we will use the invertibility of the model operators
$\mathcal L_g^{\mathcal G}$, $g\in\mathcal K_1$, with linear
transformations as well as Lemmas~\ref{lOmega'H}--\ref{lW-W}.
Then, in subsection~2 of this section, using the results
of~\cite{SkDu91} and Lemma~5.2~\cite{KovSk}, we will obtain a
priori estimates for functions with the support in the whole of
$\bar G$.

We denote ${\cal O}_\varepsilon({\cal K}_1)=\{x\in {\mathbb R}^n:
{\rm dist}(x,\ {\cal K}_1)<\varepsilon\}$.
\begin{lemma}\label{lAprLMain}
Let Conditions~\ref{condEllipPinG}--\ref{condCoordTrans} hold and,
for each $g\in{\cal K}_1$, the operator $\mathcal L_g^{\mathcal
G}$ be an isomorphism.\footnote{In subsection~5
of~\S~\ref{sectStatement}, one can find necessary and sufficient
condition under which $\mathcal L_g^{\mathcal G}$ is an
isomorphism.} Then there is an $\varepsilon$, $0<\varepsilon<{\rm
dist}({\cal K}_1,\ {\cal K}_2\cup{\cal K}_3)/2$, such that for all
$u\in\{u\in H_b^{l+2m}(G): \supp u\subset \bar G\cap {\cal
O}_\varepsilon({\cal K}_1)\}$ the following estimate holds:
 $$
   \|u\|_{H_b^{l+2m}(G)}\le c(\|{\bf L}u\|_{{\cal H}_b^l(G,\ \Upsilon)}+\|u\|_{H^0_{b+1-l-2m}(G)}),
 $$
where  $c>0$ is independent of $u$.
\end{lemma}
Using the unity partition method, Leibniz' formula,
Lemma~2.1~\cite{KondrTMMO67}, and Lemma~1.2~\cite{SkMs86}, one can
reduce the proof of Lemma~\ref{lAprLMain} to the proof of the
following result.

\begin{lemma}\label{lAprL}
Let the conditions of Lemma~\ref{lAprLMain} hold. Then for each
$g\in{\cal K}_1$ there is an
 $\varepsilon_0=\varepsilon_0(g)>0$ such that for any
$U\in\{U\in H_b^{l+2m,\,N}(\Omega): \supp U_j\subset \bar
\Omega_j\cap {\cal V}_{\varepsilon_0}(0),\
 j=1,\ \dots,\ N=N(g)\}$ the following inequality holds:
 $$
   \|U\|_{H_b^{l+2m,\,N}(\Omega)}\le c\|{\cal L}_g^\omega U\|_{H_b^{l,\,N}(\Omega)},
 $$
where ${\cal V}_{\varepsilon_0}(0)=\{x\in {\mathbb R}^n:
|x|<\varepsilon_0\},$ $c>0$ is independent of $U$.
\end{lemma}
\begin{proof}
Using the invertibility of ${\cal L}^{\cal G}_g$ and
Lemma~\ref{lBOmegaG}, for all $U\in\ H_b^{l+2m,\,N}(\Omega)$ with
$\supp U_j\subset \bar \Omega_j\cap {\cal V}_{\varepsilon_0}(0)$
we get
\begin{multline}\label{eqAprL1}
   \|U\|_{H_b^{l+2m,\,N}(\Omega)}\le k_1\|{\cal L}^{\cal G}_g U\|_{H_b^{l,\,N}(\Omega)}\le
  k_2(\|{\cal L}_g^\omega U\|_{H_b^{l,\,N}(\Omega)}+\\
 +\varepsilon_0\|U\|_{H_b^{l+2m,\,N}(\Omega)}+\sum\limits_{k=1}^N\sum\limits_{q=2}^{R_k}
 \|\zeta_{qk,3}U_k-\zeta_{qk,3}\hat U_k\|_{H_b^{l+2m}(\Omega_k)}).
\end{multline}
Let us estimate the last norm in~(\ref{eqAprL1}). By
Theorem~4.1~\cite{MP}, we have
\begin{multline}\label{eqAprL2}
  \|\zeta_{qk,3}U_k-\zeta_{qk,3}\hat U_k\|_{H_b^{l+2m}(\Omega_k)}\le k_3
 (\|{\cal P}_k(\zeta_{kq,3}U_k-\zeta_{kq,3}\hat U_k)\|_{H_b^{l}(\Omega_k)}+\\
+\|\zeta_{kq,3}U_k-\zeta_{kq,3}\hat
U_k\|_{H_{b-l-2m}^{0}(\Omega_k)}).
\end{multline}
From Lemma~\ref{lW-W} and the continuity of the embedding
$H_{b}^{l+2m}(\Omega_k)\subset H_{b-l-2m+1}^{1}(\Omega_k)$, it
follow that
\begin{equation}\label{eqAprL3}
 \|\zeta_{kq,3}U_k-\zeta_{kq,3}\hat U_k\|_{H_{b-l-2m}^{0}(\Omega_k)}\le k_4\varepsilon_0 \|U_k\|_{H_{b}^{l+2m}(\Omega_k)}.
\end{equation}
To estimate the first norm on the right-hand side
of~(\ref{eqAprL2}), we apply Leibniz' formula and
Lemmas~\ref{lOmega'H} and~\ref{lGammaOmega'}:
\begin{multline}\label{eqAprL4}
 \|{\cal P}_k(\zeta_{kq,3}U_k-\zeta_{kq,3}\hat U_k)\|_{H_b^{l}(\Omega_k)}\le
 k_5(\|\zeta_{kq,3}{\cal P}_kU_k\|_{H_b^{l}(\Omega_k)}
+\|\zeta_{kq,3}{\cal P}_k\hat U_k\|_{H_b^{l}(\Omega_k)}+\\
+\sum\limits_{|\beta|\le 2m-1}\sum\limits_{|\gamma|=2m-|\beta|}
\|D^{\gamma}\zeta_{kq,3}D^\beta U_k- D^{\gamma}\zeta_{kq,3}D^\beta\hat
U_k\|_{H_b^{l}(\Omega_k)})
\le k_6(\|{\cal P}_kU_k\|_{H_b^{l}(\Omega_k)}+\\
+\varepsilon_0\|U_k\|_{H_b^{l+2m}(\Omega_k)} + \sum\limits_{|\beta|\le
2m-1}\sum\limits_{|\gamma|=2m-|\beta|} \|D^{\gamma}\zeta_{kq,3}D^\beta U_k-
D^{\gamma}\zeta_{kq,3}D^\beta\hat U_k\|_{H_b^{l}(\Omega_k)}).
\end{multline}
Since $|D^{\gamma}\zeta_{kq,3}|\le k_7
r^{-|\gamma|}|\zeta_{kq,2}|,$ it follows that
\begin{multline}\label{eqAprL5}
 \sum\limits_{|\beta|\le 2m-1}\sum\limits_{|\gamma|=2m-|\beta|}
\|D^{\gamma}\zeta_{kq,3}D^\beta U_k- D^{\gamma}\zeta_{kq,3}D^\beta\hat U_k\|_{H_b^{l}(\Omega_k)}\le\\
  \le k_8\sum\limits_{|\alpha|\le l+2m-1}
\|\zeta_{kq,2}D^{\alpha} U_k-\zeta_{kq,2} D^{\alpha} \hat U_k\|_{H_{b+|\alpha|-l-2m}^{0}(\Omega_k)}\le\\
\le k_{9}\sum\limits_{|\alpha|\le l+2m-1}
\{\|\zeta_{kq,2}D^{\alpha}U_k-\zeta_{kq,2}\widehat{D^{\alpha}U_k}\|_{H_{b+|\alpha|-l-2m}^{0}(\Omega_k)}+\\
+\|\zeta_{kq,2}\widehat{D^{\alpha} U_k}-\zeta_{kq,2}D^{\alpha}\hat
U_k\|_{H_{b+|\alpha|-l-2m}^{0}(\Omega_k)} \}.
\end{multline}
Using Lemma~\ref{lW-W} and the continuity of the embedding
$H_b^{l+2m}(\Omega_k)\subset
H_{b+1+|\alpha|-l-2m}^{1+|\alpha|}(\Omega_k)$ for $|\alpha|\le
l+2m-1$, we obtain
\begin{multline}\label{eqAprL6}
 \|\zeta_{kq,2}D^{\alpha}U_k-\zeta_{kq,2}\widehat{D^{\alpha}U_k}\|_{H_{b+|\alpha|-l-2m}^{0}(\Omega_k)}\le\\
\le
k_{10}\varepsilon_0\|D^{\alpha}U_k\|_{H_{b+1+|\alpha|-l-2m}^{1}(\Omega_k)}
 \le k_{11}\varepsilon_0\|U_k\|_{H_{b}^{l+2m}(\Omega_k)}.
\end{multline}
Similarly, from Lemma~\ref{lGammaOmega'}, it follows that
\begin{equation}\label{eqAprL7}
\|\zeta_{kq,2}\widehat{D^{\alpha} U_k}-\zeta_{kq,2}D^{\alpha}\hat
U_k\|_{H_{b+|\alpha|-l-2m}^{0}(\Omega_k)}\le
k_{12}\varepsilon_0\|U_k\|_{H_b^{l+2m}(\Omega_k)}.
\end{equation}
Now the conclusion of the lemma follows
from~(\ref{eqAprL1})--(\ref{eqAprL7}) with sufficiently small
$\varepsilon_0$.
\end{proof}

\medskip

{\bf 2.} Repeating the proof of Theorem~2.1~\cite{SkDu91} and
taking into account Lemma~5.2~\cite{KovSk}, from
Lemma~\ref{lAprLMain} of the present work and Lemmas~2.4 and~2.5
of~\cite{SkDu91}, we deduce the following result.
\begin{theorem}\label{thAprLMain}
Let the conditions of Lemma~\ref{lAprLMain} hold and $b>l+2m-1$.
Then, for all $u\in H_b^{l+2m}(G)$, the following estimate holds:
\begin{equation}\label{eqAprLMain}
  \|u\|_{H_b^{l+2m}(G)}\le c(\|{\bf L}u\|_{{\cal H}_b^l(G,\ \Upsilon)}+\|u\|_{H^0_{b+1-l-2m}(G)}),
\end{equation}
where $c>0$ is independent of $u$.
\end{theorem}

By virtue of the compactness of the embedding
$H_b^{l+2m}(G)\subset H^0_{b+1-l-2m}(G)$ (see
Lemma~3.5~\cite{KondrTMMO67}), from Theorem~\ref{thAprLMain} it
follows that the operator ${\bf L}$ has a finite-dimensional
kernel and a closed range.

%% file: sect5.tex
\section{Construction of right regularizer}\label{sectRightReg}
In this section, we construct a right regularizer for ${\bf L}$,
which, being combined with Theorem~\ref{thAprLMain}, allows us to
prove the Fredholm solvability of nonlocal boundary-value
problem~(\ref{eqPinG}), (\ref{eqBinG}).

{\bf 1.} To begin with, we consider the case where the supports of
functions are subsets of a neighborhood of ${\cal K}_1$. In this
situation, we will use the invertibility of the operators
$\mathcal L_g^{\mathcal G}$, $g\in\mathcal K_1$, with linear
transformations as well as some special constructions
``compensating'' the nonlinearity in the argument transformations.
Then, in subsection~2 of this section, using the results
of~\cite{SkDu91} and Lemma~5.2~\cite{KovSk}, we will construct the
right regularizer in the whole of $G$.

First of all, let us prove the following auxiliary result.
\begin{lemma}\label{lSmallComp}
Let $H,\ H_1$, and $H_2$ be Hilbert spaces, ${\cal A}: H\to H_1$ a
linear bounded operator, ${\cal T}_0: H\to H_2$ a linear compact
operator. Suppose that, for some $\varepsilon,\ c>0$ and all $f\in
H$, the following inequality holds:
 \begin{equation}\label{eqSmallComp}
  \|{\cal A}f\|_{H_1}\le \varepsilon\|f\|_H+c\|{\cal T}_0 f\|_{H_2}.
 \end{equation}
Then there are bounded operators ${\cal M},\ {\cal F}:H\to H_1$
such that
 $$
  {\cal A}={\cal M}+{\cal F},
 $$
where $\|{\cal M}\|\le 2\varepsilon$ and the operator ${\cal F}$
is finite-dimensional.
\end{lemma}
\begin{proof}
As is well known (see, e.g.,~\cite[Chapter~5, Section~85]{RS}),
any compact operator is the limit of a uniformly convergent
sequence of finite-dimensional operators. Therefore, there are
bounded operators ${\cal M}_0,\ {\cal F}_0:H\to H_2$ such that
${\cal T}_0={\cal M}_0+{\cal F}_0$, $\|{\cal M}_0\|\le
c^{-1}\varepsilon$, and ${\cal F}_0$ is finite-dimensional. From
this and~(\ref{eqSmallComp}), it follows that
\begin{equation}\label{eqSmallComp1}
\|{\cal A}f\|_{H_1}\le 2\varepsilon\|f\|_H+c\|{\cal F}_0
f\|_{H_2}\quad \mbox{for all } f\in H.
\end{equation}
We denote by $\ker({\cal F}_0)^\bot$ the orthogonal supplement in
$H$ to the kernel of ${\cal F}_0$. Since the finite-dimensional
operator ${\cal F}_0$ maps $\ker({\cal F}_0)^\bot$ onto its range
in a one-to-one manner, it follows that the subspace $\ker({\cal
F}_0)^\bot$ is of finite dimension. Let ${\cal I}$ denote the
unity operator in $H$ and ${\cal P}_0$ the orthogonal projector
onto $\ker({\cal F}_0)^\bot$. Obviously, ${\cal A}{\cal P}_0:H\to
H_1$ is a finite-dimensional operator. Furthermore, since ${\cal
I}-{\cal P}_0$ is the orthogonal projector onto $\ker({\cal
F}_0)$, it follows that ${\cal F}_0({\cal I}-{\cal P}_0)=0$.
Therefore, substituting in~(\ref{eqSmallComp1}) the function
$({\cal I}-{\cal P}_0)f$ for $f$, we get
$$
\|{\cal A}({\cal I}-{\cal P}_0) f\|_{H_1}\le 2\varepsilon\|({\cal
I}-{\cal P}_0)f\|_H\le2\varepsilon\|f\|_H\quad \text{for all }
  f\in H.
$$
Denoting ${\cal M}={\cal A}({\cal I}-{\cal P}_0)$ and ${\cal
F}={\cal A}{\cal P}_0$ completes the proof.
\end{proof}

Now we proceed to construct the right regularizer.

\begin{lemma}\label{lRegGK1}
Let the conditions of Lemma~\ref{lAprLMain} hold. Then, for all
sufficiently small $\varepsilon$, $0<\varepsilon<{\rm dist}({\cal
K}_1,\ {\cal K}_2\cup{\cal K}_3)/2$, there are bounded operators
${\bf R}_1$, ${\bf M}_1$ and a compact operator ${\bf T}_1$ acting
from $\{f\in {\cal H}_b^l(G,\ \Upsilon): \supp f\subset \bar G\cap
{\cal O}_\varepsilon({\cal K}_1)\}$ to $H_b^{l+2m}(G)$, ${\cal
H}_b^l(G,\ \Upsilon)$, and ${\cal H}_b^l(G,\ \Upsilon)$
respectively and such that
 $$
  {\bf L}{\bf R}_1 f=f+{\bf M}_1f+{\bf T}_1f,
 $$
$\|{\bf M}_1f\|_{{\cal H}_b^l(G,\ \Upsilon)}\le
c\varepsilon\|f\|_{{\cal H}_b^l(G,\ \Upsilon)}$. Here $c>0$ is
independent of $\varepsilon$ and $f$.
\end{lemma}

Using the unity partition method, Leibniz' formula, and
Lemma~2.1~\cite{KondrTMMO67}, one can reduce the proof of
Lemma~\ref{lRegGK1} to the proof of the following result.

\begin{lemma}\label{lRegOmega}
Let the conditions of Lemma~\ref{lAprLMain} hold. Then, for each
$g\in{\cal K}_1$ and all sufficiently small
$\varepsilon_1=\varepsilon_1(g)>0$, there are bounded operators
${\cal R}_g$,  ${\cal M}_g$ and a compact operator ${\cal T}_g$
acting from $\{f\in {\cal H}_b^{l,\,N}(\Omega,\ \Gamma):
 \supp f\subset {\cal V}_{\varepsilon_1}(0)\}$
to $H_b^{l+2m,\,N}(\Omega)$, ${\cal H}_b^{l,\,N}(\Omega,\ \Gamma)$
and ${\cal H}_b^{l,\,N}(\Omega,\ \Gamma)$ respectively and such
that
 \begin{equation}\label{eqRegOmega}
  {\cal L}_g^\omega{\cal R}_g f=f+{\cal M}_g f+
  {\cal T}_g f,
 \end{equation}
$\|{\cal M}_g f\|_{H_b^l(G,\ \Gamma)}\le
c\varepsilon_1\|f\|_{H_b^l(G,\ \Gamma)}$. Here $c>0$ is
independent of $\varepsilon_1$ and $f$.
\end{lemma}
\begin{proof}
1) As before, we denote $d_1=\min\{1,\ \chi_{j\sigma kqs}\}/2,$
$d_2=2\max\{1,\ \chi_{j\sigma kqs}\}$. We choose
$\varepsilon_1<d_1\varepsilon_0/4,$ where $\varepsilon_0$ is
defined in Lemma~\ref{lAprL}. We introduce a function
$\psi_{\varepsilon_1}(x)=\psi(x/\varepsilon_1),$ where $\psi\in
C^\infty({\mathbb R}^n)$, $\psi(x)=1$ for $|x|\le1$, $\psi(x)=0$
for $|x|\ge 2$. It is obvious that $\psi_{\varepsilon_1}\in
C^\infty({\mathbb R}^n)$, $\psi_{\varepsilon_1}(x)=1$ for $|x|\le
\varepsilon_1$, $\psi_{\varepsilon_1}(x)=0$ for $|x|\ge
2\varepsilon_1$. Since $|D^\alpha \psi_{\varepsilon_1}|\le
c_{\alpha}r^{-|\alpha|},$ from Lemma~2.1~\cite{KondrTMMO67} it
follows that
\begin{equation}\label{eqPsi}
\|\psi_{\varepsilon_1}v\|_{H_b^{l+2m}(\Omega_k)}\le c
\|v\|_{H_b^{l+2m}(\Omega_k)}\quad \mbox{for all }
   v\in H_b^{l+2m}(\Omega_k),
\end{equation}
where $c>0$ is independent of $\varepsilon_1$. Moreover, we assume
that $\psi_{\varepsilon_1}$, being written in cylindrical
coordinates, does not depend on $\varphi$.

Put $f_0=\{f_j\}$, $g=\{g_{j\sigma\mu}\}$, $\{f_0,\ g\}=\{f_j,\
g_{j\sigma\mu}\}$.

By assumption, the operator ${\cal L}_g^{\cal G}:
H_b^{l+2m,\,N}(\Omega)\to {\cal H}_b^{l,\,N}(\Omega,\ \Gamma)$ has
a bounded inverse $({\cal L}_g^{\cal G})^{-1}: {\cal
H}_b^{l,\,N}(\Omega,\ \Gamma)\to H_b^{l+2m,\,N}(\Omega).$
Therefore, we can introduce the operators
$$
{\cal R}_1:H_b^{l,\,N}(\Omega)\to H_b^{l+2m,\,N}(\Omega),\quad
{\cal R}_2:{\cal H}_b^{l,\,N}(\Gamma)\to H_b^{l+2m,\,N}(\Omega)
$$
given by
$$
 {\cal R}_1f_0=\psi_{\varepsilon_1}({\cal L}_g^{\cal G})^{-1}\{f_0,\
 0\},\quad
 {\cal R}_2g=\psi_{\varepsilon_1}({\cal L}_g^{\cal G})^{-1}\{0,\
 g\},
$$
where ${\cal H}_b^{l,\,N}(\Gamma)= \prod\limits_{j,\,\sigma,\,\mu}
H_b^{l+2m-m_{j\sigma\mu}-1/2}(\Gamma_{j\sigma})$. Thus, the
supports of ${\cal R}_1f_0$ and ${\cal R}_2g$ are subsets of the
ball of radius~$2\varepsilon_1$ centered at the origin.

Let us introduce the operators
\begin{align*}
{\cal P}:\ &H_b^{l+2m,\,N}(\Omega)\to H_b^{l,\,N}(\Omega),\\
 {\cal B}^{\cal G},\ {\cal B}^{\omega}:\
&H_b^{l+2m,\,N}(\Omega)\to {\cal H}_b^{l,\,N}(\Gamma)
\end{align*}
given by
$$
 {\cal P}U=\{{\cal P}_jU_j\},\quad
  {\cal B}^{\cal G}U=\{{\cal B}^{\cal G}_{j\sigma\mu}U\},\quad
{\cal B}^{\omega}U=\{{\cal B}^{\omega}_{j\sigma\mu}U\}.
$$

Now we establish a relation between the operators ${\cal P}$,
${\cal B}^{\cal G}$, ${\cal B}^\omega$ and ${\cal R}_1$, ${\cal
R}_2$. To this end, we will use the following well-known property
of weighted spaces (see Lemma~3.5~\cite{KondrTMMO67}): $(*)$~{\it
the embedding operator from $\{v\in H_b^{l+1}(\Omega_j):\supp
v\subset{\cal V}_d(0),\ d>0\}$ into $H_b^{l}(\Omega_j)$ is
compact.}

From Leibniz' formula, the boundedness of
$\supp\psi_{\varepsilon_1}$, and property~$(*)$, it follows that
\begin{equation}\label{eqRegOmega1}
 {\cal P}{\cal R}_1f_0=\psi_{\varepsilon_1} f_0+{\cal
 T}_{1}f_0,\quad
 {\cal P}{\cal R}_2g={\cal T}_{2}g,
\end{equation}
where ${\cal T}_{1}:H_b^{l,\,N}(\Omega)\to H_b^{l,\,N}(\Omega)$
and ${\cal T}_{2}:{\cal H}_b^{l,\,N}(\Gamma) \to
H_b^{l,\,N}(\Omega)$ are compact operators. Similarly,
\begin{multline}\label{eqRegOmega2}
{\cal B}^{\cal G}{\cal R}_2g=\psi_{\varepsilon_1} g+\\
+\Big\{\sum\limits_{k,\,q,\,s}( \psi_{\varepsilon_1}(\chi_{j\sigma
kqs}x)-\psi_{\varepsilon_1}(x) )
 (B_{j\sigma\mu kqs}[({\cal L}_g^{\cal G})^{-1}\{0,\ g\}]_k)
({\cal G}_{j\sigma kqs}y,\ z)|_{\Gamma_{j\sigma}}\Big\}+{\cal
T}_{3}g,
\end{multline}
where ${\cal T}_{3}$ is a compact operator in ${\cal
H}_b^{l,\,N}(\Gamma)$; here and in what follows, we denote by
$[\cdot]_{k}$ the $k$th component of an $N$-dimensional vector and
by $\{\dots\}$ a vector with the components defined by the indices
$j,\ \sigma,\ \mu$.

Let us show that each term in the sum in~(\ref{eqRegOmega2}) is a
compact operator. Let $\zeta_{kq,i}$ be the functions defined by
formulas~(\ref{eqZetaJSigma}). We also introduce the functions
 $\hat\psi_0,$ $\hat\psi_1\in C_0^\infty({\mathbb R}^n)$ such that
$$
 \hat\psi_1(x)=1\ \mbox{for } 2d_1\varepsilon_1\le |x|\le d_2 \varepsilon_1,\
 \hat\psi_1(x)=0\ \mbox{outside } d_1\varepsilon_1\le |x|\le 2d_2 \varepsilon_1,
$$
$$
 \hat\psi_0(x)=1\ \mbox{for } d_1\varepsilon_1\le |x|\le 2 d_2 \varepsilon_1,\
 \hat\psi_0(x)=0\ \mbox{outside } d_1\varepsilon_1/2\le |x|\le 4 d_2 \varepsilon_1.
$$
Then, by virtue of the boundedness of the trace operator in
weighted spaces, we have
\begin{multline}\label{eqRegOmega3}
 \|( \psi_{\varepsilon_1}(\chi_{j\sigma kqs}x)-\psi_{\varepsilon_1}(x) )
(B_{j\sigma\mu kqs}[({\cal L}_g^{\cal G})^{-1}\{0,\ g\}]_k)
({\cal G}_{j\sigma kqs}y,\ z)|_{\Gamma_{j\sigma}}\|_{H_b^{l+2m-m_{j\sigma\mu}-1/2}(\Gamma_{j\sigma})}\le\\
\le k_2\|\zeta_{kq,2}(
\psi_{\varepsilon_1}(x)-\psi_{\varepsilon_1}(\chi_{j\sigma
kqs}^{-1}x) ) B_{j\sigma\mu kqs}[({\cal L}_g^{\cal G})^{-1}\{0,\
g\}]_k
\|_{H_b^{l+2m-m_{j\sigma\mu}}(\Omega_k)}\le\\
\le k_3\|\zeta_{kq,1}\hat\psi_1 [({\cal L}_g^{\cal G})^{-1}\{0,\
g\}]_k \|_{H_b^{l+2m}(\Omega_k)}.
\end{multline}
Since the support of $\hat\psi_1$ is bounded and does not
intersect with the origin and $\zeta_{kq,1}$ vanishes near the
sides of the angle $\Omega_k$, we can apply
Theorem~5.1~\cite[Chapter~2]{LM}. Then, using the relation ${\cal
P}_k[({\cal L}_g^{\cal G})^{-1}\{0,\ g\}]_k=0$,
from~(\ref{eqRegOmega3}) we get
 \begin{multline}\notag
 \|( \psi_{\varepsilon_1}(\chi_{j\sigma kqs}x)-\psi_{\varepsilon_1}(x) )(B_{j\sigma\mu kqs}[({\cal L}_g^{\cal G})^{-1}\{0,\ g\}]_k)
({\cal G}_{j\sigma kqs}y,\ z)|_{\Gamma_{j\sigma}}\|_{H_b^{l+2m-m_{j\sigma\mu}-1/2}(\Gamma_{j\sigma})}\le\\
\le k_4\|\hat\psi_0 [({\cal L}_g^{\cal G})^{-1}\{0,\
g\}]_k\|_{H_b^{l+2m-1}(\Omega_k)}.
 \end{multline}
Since the support of $\hat\psi_0$ is bounded, from the last
inequality and property~$(*)$ it follows that
$$
 \Big\{\sum\limits_{k,\,q,\,s}( \psi_{\varepsilon_1}(\chi_{j\sigma kqs}x)-\psi_{\varepsilon_1}(x) )
(B_{j\sigma\mu kqs}[({\cal L}_g^{\cal G})^{-1}\{0,\ g\}]_k) ({\cal
G}_{j\sigma kqs}y,\ z)|_{\Gamma_{j\sigma}}\Big\}
$$
is a compact operator acting in ${\cal H}_b^{l,\,N}(\Gamma)$.
Combining this with~(\ref{eqRegOmega2}) yields
\begin{equation}\label{eqRegOmega4}
{\cal B}^{\cal G}{\cal R}_2g=\psi_{\varepsilon_1} g+{\cal T}_{4}g,
\end{equation}
where ${\cal T}_{4}$ is a compact operator acting in ${\cal
H}_b^{l,\,N}(\Gamma)$.

Finally, from~(\ref{eqRegOmega4}), we obtain the formula for the
composition ${\cal B}^\omega{\cal R}_2$:
\begin{multline}\label{eqRegOmega5}
{\cal B}^\omega{\cal R}_2g=\psi_{\varepsilon_1} g+{\cal T}_{4}g+\\
+\Big\{\sum\limits_{k,\,q,\,s}\Bigl( (B_{j\sigma\mu kqs}[{\cal
R}_2g]_k)(\omega'_{j\sigma kqs}(y,\ z),\ z)|_{\Gamma_{j\sigma}}-
(B_{j\sigma\mu kqs}[{\cal R}_2g]_k)({\cal G}_{j\sigma kqs}y,\
z)|_{\Gamma_{j\sigma}}\Bigr)\Big\}.
\end{multline}

2) Let us introduce the operator ${\cal R}_g:{\cal
H}_b^{l,\,N}(\Omega,\ \Gamma)\to H_b^{l+2m,\,N}(\Omega)$ given by
$$
 {\cal R}_g\{f_0,\ g\}={\cal R}_1f_0-{\cal R}'_2{\cal B}^\omega {\cal R}_1f_0+{\cal R}_2g.
$$
Here ${\cal R}'_2:{\cal H}_b^{l,\,N}(\Gamma)\to
H_b^{l+2m,\,N}(\Omega)$ is the bounded operator given by
$$
 {\cal R}'_2g=\psi_{\varepsilon_1}(d_1x/2)({\cal L}_g^{\cal G})^{-1}\{0,\ g\}.
$$
Similarly to~(\ref{eqRegOmega1}) and (\ref{eqRegOmega5}), one can
prove that
\begin{equation}\label{eqRegOmega6}
{\cal P}{\cal R}'_2g={\cal T}'_{2}g,
\end{equation}
\begin{multline}\label{eqRegOmega7}
{\cal B}^\omega{\cal R}'_2g=\psi_{\varepsilon_1}(d_1x/2)g+
{\cal T}'_{4}g+\\
+\Big\{\sum\limits_{k,\,q,\,s}\Bigl( (B_{j\sigma\mu kqs}[{\cal
R}'_2g]_k)(\omega'_{j\sigma kqs}(y,\ z),\ z)|_{\Gamma_{j\sigma}}-
(B_{j\sigma\mu kqs}[{\cal R}'_2g]_k)({\cal G}_{j\sigma kqs}y,\
z)|_{\Gamma_{j\sigma}}\Bigr)\Big\},
\end{multline}
where ${\cal T}'_{2},$ ${\cal T}'_{4}$ are compact operators
acting in the same spaces as the operators ${\cal T}_{2},$ ${\cal
T}_{4}$ do.

Let us show that the operator ${\cal R}_g$ satisfies
relation~(\ref{eqRegOmega}). From~(\ref{eqRegOmega1}) and
(\ref{eqRegOmega6}), it follows that
\begin{equation}\label{eqRegOmega8}
 {\cal P}{\cal R}_g\{f_0,\ g\}=\psi_{\varepsilon_1} f_0+{\cal T}_{5}\{f_0,\ g\},
\end{equation}
where ${\cal T}_{5}:{\cal H}_b^{l,\,N}(\Omega,\ \Gamma)\to
H_b^{l,\,N}(\Omega)$ is a compact operator.

Taking into account that $\psi_{\varepsilon_1}(d_1x/2){\cal
B}^\omega{\cal R}_1f_0\equiv {\cal B}^\omega{\cal R}_1f_0$ and
using~(\ref{eqRegOmega7}), we derive
 \begin{multline}\notag
{\cal B}^\omega {\cal R}_g\{f_0,\ g\}={\cal B}^\omega{\cal
R}_1f_0- {\cal B}^\omega{\cal R}'_2{\cal B}^\omega {\cal R}_1f_0+
{\cal B}^\omega {\cal R}_2g=\\
=-{\cal T}'_{4}{\cal B}^\omega {\cal R}_1f_0-
\Big\{\sum\limits_{k,\,q,\,s}\Bigl( (B_{j\sigma\mu kqs}[{\cal
R}'_2{\cal B}^\omega
         {\cal R}_1f_0]_k)(\omega'_{j\sigma kqs}(y,\ z),\ z)|_{\Gamma_{j\sigma}}-\\
-(B_{j\sigma\mu kqs}[{\cal R}'_2{\cal B}^\omega {\cal
R}_1f_0]_k)({\cal G}_{j\sigma kqs}y,\
z)|_{\Gamma_{j\sigma}}\Bigr)\Big\}+ {\cal B}^\omega {\cal R}_2g.
 \end{multline}
From this, using~(\ref{eqRegOmega5}), we obtain
\begin{multline}\label{eqRegOmega9}
 {\cal B}^\omega {\cal R}_g g=\psi_{\varepsilon_1} g+
{\cal T}_{6}\{f_0,\ g\}+\\
+\Big\{\sum\limits_{k,\,q,\,s}\Bigl( (B_{j\sigma\mu kqs}[{\cal
R}_2g]_k)(\omega'_{j\sigma kqs}(y,\ z),\ z)|_{\Gamma_{j\sigma}}
-(B_{j\sigma\mu kqs}[{\cal R}_2g]_k)({\cal G}_{j\sigma kqs}y,\ z)|_{\Gamma_{j\sigma}}\Bigr)\Big\}-\\
-\Big\{\sum\limits_{k,\,q,\,s}\Bigl( (B_{j\sigma\mu kqs}[{\cal
R}'_2{\cal B}^\omega
         {\cal R}_1f_0]_k)(\omega'_{j\sigma kqs}(y,\ z),\ z)|_{\Gamma_{j\sigma}}-\\
-(B_{j\sigma\mu kqs}[{\cal R}'_2{\cal B}^\omega {\cal
R}_1f_0]_k)({\cal G}_{j\sigma kqs}y,\
z)|_{\Gamma_{j\sigma}}\Bigr)\Big\},
\end{multline}
where ${\cal T}_{6}:{\cal H}_b^{l,\,N}(\Omega,\ \Gamma)\to {\cal
H}_b^{l,\,N}(\Gamma)$ is a compact operator.

Let us consider the terms of the first sum on the right-hand side
of~(\ref{eqRegOmega9}). By Lemma~\ref{lBOmegaG}, we have
 \begin{multline}\label{eqRegOmega10}
  \|(B_{j\sigma\mu kqs}[{\cal R}_2g]_k)(\omega'_{j\sigma kqs}(y,\ z),\ z)|_{\Gamma_{j\sigma}}-\\
-(B_{j\sigma\mu kqs}[{\cal R}_2g]_k)({\cal G}_{j\sigma kqs}y,\
z)|_{\Gamma_{j\sigma}}\|_
      {H_b^{l+2m-m_{j\sigma\mu}-1/2}(\Gamma_{j\sigma})}\le
 k_5(\varepsilon_1\|[{\cal R}_2g]_k\|_{H_b^{l+2m}(\Omega_k)}+\\
 +\|\zeta_{kq,3}[{\cal R}_2g]_k-
  \zeta_{kq,3}\widehat{[{\cal R}_2g]_k}\|_{H_b^{l+2m}(\Omega_k)}).
 \end{multline}

From inequalities~(\ref{eqAprL2})--(\ref{eqAprL7}) for the
function $U_k=[{\cal R}_2g]_k$, inequality~(\ref{eqRegOmega10}),
and the second relation in~(\ref{eqRegOmega1}), we obtain
  \begin{multline}\notag
  \|(B_{j\sigma\mu kqs}[{\cal R}_2g]_k)(\omega'_{j\sigma kqs}(y,\ z),\ z)|_{\Gamma_{j\sigma}}-\\
-(B_{j\sigma\mu kqs}[{\cal R}_2g]_k)({\cal G}_{j\sigma kqs}y,\
z)|_{\Gamma_{j\sigma}}\|_
      {H_b^{l+2m-m_{j\sigma\mu}-1/2}(\Gamma_{j\sigma})}\le\\
\le k_6(\varepsilon_1\|[{\cal R}_2g]_k\|_{H_b^{l+2m}(\Omega_k)}+\|{\cal P}_k[{\cal R}_2g]_k\|_{H_b^{l}(\Omega_k)})=\\
=k_6(\varepsilon_1\|\psi_{\varepsilon_1}[({\cal L}_g^{\cal
G})^{-1}\{0,\ g\}]_k\|_{H_b^{l+2m}(\Omega_k)}+ \|[{\cal
T}_2g]_k\|_{H_b^{l}(\Omega_k)}).
 \end{multline}
This, being combined with inequality~(\ref{eqPsi}) and the
boundedness of the operator $({\cal L}_g^{\cal G})^{-1}: {\cal
H}_b^{l,\,N}(\Omega,\ \Gamma)\to H_b^{l+2m,\,N}(\Omega)$, finally
implies
  \begin{multline}
  \|(B_{j\sigma\mu kqs}[{\cal R}_2g]_k)(\omega'_{j\sigma kqs}(y,\ z),\ z)|_{\Gamma_{j\sigma}}-\\
-(B_{j\sigma\mu kqs}[{\cal R}_2g]_k)({\cal G}_{j\sigma kqs}y,\
z)|_{\Gamma_{j\sigma}}\|_
      {H_b^{l+2m-m_{j\sigma\mu}-1/2}(\Gamma_{j\sigma})}\le\\
\le k_7 (\varepsilon_1\|g\|_{{\cal H}_b^{l,\,N}(\Gamma)} +\|[{\cal
T}_2g]_k\|_{H_b^{l}(\Omega_k)}).
 \end{multline}
Therefore, by Lemma~\ref{lSmallComp}, we have
 \begin{multline}\notag
 (B_{j\sigma\mu kqs}[{\cal R}_2g]_k)(\omega'_{j\sigma kqs}(y,\ z),\ z)|_{\Gamma_{j\sigma}}
-(B_{j\sigma\mu kqs}[{\cal R}_2g]_k)({\cal G}_{j\sigma kqs}y,\ z)|_{\Gamma_{j\sigma}}=\\
={\cal M}_{j\sigma\mu kqs}g+{\cal F}_{j\sigma\mu kqs}g
 \end{multline}
with the operators
$$
 {\cal M}_{j\sigma\mu kqs},\ {\cal F}_{j\sigma\mu kqs}:
{\cal H}_b^{l,\,N}(\Gamma)\to
H_b^{l+2m-m_{j\sigma\mu}-1/2}(\Gamma_{j\sigma})
$$
such that $\|{\cal M}_{j\sigma\mu kqs}\|\le 2k_7\varepsilon_1$ and
the operator ${\cal F}_{j\sigma\mu kqs}$ is finite-dimensional.

Analogously, one can prove that each term of the second sum on the
right-hand side of~(\ref{eqRegOmega9}) can be represented as the
sum of an operator with small norm and a compact one. From this,
(\ref{eqRegOmega9}), and~(\ref{eqRegOmega8}), choosing
$\supp\{f_0,\ g\}\subset{\cal V}_{\varepsilon_1}(0)$, we get the
conclusion of the lemma.
\end{proof}

\medskip

{\bf 2.} Now we can prove that, under certain conditions, the
operator ${\bf L}:H_b^{l+2m}(G)\to {\cal H}_b^l(G,\ \Upsilon)$ is
Fredholm.

\begin{theorem}\label{thFredLinG}
Let the conditions of Lemma~\ref{lAprLMain} hold and $b>l+2m-1$.
Then the operator ${\bf L}:H_b^{l+2m}(G)\to  {\cal H}_b^l(G,\
\Upsilon)$ is Fredholm.
\end{theorem}
\begin{proof}
By virtue of Theorem~\ref{thAprLMain} of the present paper and
Theorems~7.1, 15.2~\cite{Kr}, it suffices to construct a right
regularizer ${\bf R}$ for ${\bf L}$.

Repeating the arguments of~\cite[\S~3]{SkDu91} and taking into
account Lemma~5.2~\cite{KovSk}, from Lemma~\ref{lRegGK1} of the
present paper we deduce the existence of bounded operators
\begin{align*}
 {\bf R}'&:{\cal H}_b^l(G,\ \Upsilon)\to H_b^{l+2m}(G),\\
 {\bf M},\ {\bf T}&:{\cal H}_b^l(G,\ \Upsilon)\to {\cal
H}_b^l(G,\ \Upsilon)
\end{align*}
such that
$$
 \mathbf L{\bf R}'=\mathbf I+{\bf M}+{\bf T},
$$
where $\|\mathbf M\|<1$ and the operator ${\bf T}$ is compact.
Since $\|\mathbf M\|<1$, it follows that the operator $\mathbf
I+{\bf M}$ has a bounded inverse. Obviously, the operator ${\bf
R}={\bf R}'({\bf I}+{\bf M})^{-1}$ is a right regularizer
for~${\bf L}$.
\end{proof}

\medskip

{\bf 3.} Until now, we assumed that $b>l+2m-1$. In this
subsection, using results of~\cite{SkMs86}, we study the case
where $b$ is arbitrary but $n=2$. As mentioned before, if $b$ is
arbitrary, we have to consider solutions and right-hand sides of
the nonlocal problem as functions with power singularities not
only near the set $\mathcal K_1$ but also near $\mathcal K_2$ and
$\mathcal K_3$. This corresponds to the consistency conditions
(see~\S~\ref{sectStatement}).

Thus, let $n=2$. We introduce the space $\tilde H_b^l(G)$ as the
completion of $C_0^\infty(\bar G\setminus \mathcal K)$ with
respect to the norm
$$
\|u\|_{H_b^l(G)}=\left(\sum\limits_{|\alpha|\le l}\int\limits_G
\tilde\rho^{2(b-l+|\alpha|)}|D^\alpha u|^2 dy\right)^{1/2},
$$
where $\tilde\rho=\tilde\rho(y)={\rm dist}(y,\ {\cal K})$
(cf.~\S~\ref{sectStatement}). For $l\ge1$, we denote by $\tilde
H_b^{l-1/2}(\Upsilon)$ the space of traces on a smooth curve
$\Upsilon\subset\bar G$ with the norm
$$
\|\psi\|_{\tilde H_b^{l-1/2}(\Upsilon)}=\inf\|u\|_{\tilde
H_b^l(G)}\quad (u\in \tilde H_b^l(G):\ u|_\Upsilon=\psi).
$$

We assume that the following condition holds.
\begin{condition}\label{condKinG}
If $g\in{\cal K}_{3}\cap\omega_{is}(\Upsilon_i)\ne\varnothing,$
then $\omega_{is}^{-1}(g)\in{\cal K}.$
\end{condition}

The fulfillment of Condition~\ref{condKinG} guarantees that the
set of points in which the consistency condition must be imposed
is finite. If Condition~\ref{condKinG} fails, then the consecutive
shifts of the set ${\cal K}_1$ (under the transformations
$\omega_{is}$ and $\omega_{is}^{-1}$) may form an infinite set,
which should be used instead of ${\cal K}$ in the definition of
weighted spaces.

\smallskip

In this subsection, we consider the following bounded operator
corresponding to problem\footnote{Notice that
equation~(\ref{eqPinG}) is now considered in $G\setminus\mathcal
K_3$ but not in the whole of $G$.}~(\ref{eqPinG}), (\ref{eqBinG}):
$$
 {\bf L}=\{{\bf P}(y,\ D),\ {\bf B}_{i\mu}(y,\ D)\}: \tilde H_b^{l+2m}(G)\to
 \tilde H_b^l(G)\times\prod\limits_{i=1}^{N_0}\prod\limits_{\mu=1}^m
\tilde H_b^{l+2m-m_{i\mu}-1/2}(\Upsilon_i),\quad b\in\mathbb R.
$$

Since solutions and right-hand sides of the nonlocal problem may
now have power singularities near the points of $\mathcal K_2$ and
$\mathcal K_3$, we have to consider the model problems
corresponding to these points in weighted spaces but not in the
Sobolev spaces.

We fix a point $g\in{\cal K}_2\cup{\cal K}_3.$ Let $y\mapsto
y'(g)$ be a non-degenerate infinitely differentiable argument
transformation mapping some neighborhood ${\cal V}(g)$ of the
point $g$ onto a neighborhood ${\cal V}_g(0)$ of the origin, so
that the point $g$ maps to the origin. We denote by ${\cal
P}(D_y),$ ${\cal B}_{i\mu 0}(D_y)$ the principal homogeneous parts
of the operators ${\bf P}(g,\ D),$ $B_{i\mu 0}(g,\ D)$ written in
new coordinates $y'=y'(g)$ (with after-denoting $y'$ by $y$). Now
we write the operators ${\cal P}(D_y),$ ${\cal B}_{i\mu 0}(D_y)$
in polar coordinates: ${\cal P}(D_y)=r^{-2m}\tilde{\cal
P}(\varphi,\ D_\varphi,\ rD_r),$ ${\cal
B}_{i\mu0}(D_y)=r^{-m_{i\mu}}\tilde{\cal B}_{i\mu0}(\varphi,\
D_\varphi,\ rD_r).$

If $g\in{\cal K}_{2},$ then $g\in\Upsilon_i$ for some $i=i(g)$. By
virtue of the smoothness of $\Upsilon_i$, in a sufficiently small
neighborhood ${\cal V}(g)$ of $g$ there is a non-degenerate
infinitely smooth argument transformation $y\mapsto y'=y'(g)$
mapping ${\cal V}(g)\cap G$ onto the intersection of the
half-plane ${\mathbb R}_+^2=\{y: |\varphi|<\pi/2\}$ with a
neighborhood of ${\cal V}_g(0)$. Let us introduce the bounded
operator
$$
 {\cal L}_g:H_b^{l+2m}(K_{\pi/2})\to H_b^l(K_{\pi/2})\times\prod\limits_{j=1}^2
 \prod\limits_{\mu=1}^m H_b^{l+2m-m_{i\mu}-1/2}(\gamma_j)
$$
given by
$$
 {\cal L}_gU=\{{\cal P}(D_y)U,\ {\cal B}_{i\mu 0}(D_y)U|_{\gamma_j}\},
$$
where $K_{\pi/2}=\{y: |\varphi|<\pi/2\},$ $\gamma_j=\{y:
\varphi=(-1)^j\pi/2\},$ $j=1,\ 2.$ We also introduce the bounded
operator
$$
 \tilde{\cal L}_g(\lambda):W_2^{l+2m}(-\pi/2,\ \pi/2)\to {\cal W}_2^l[-\pi/2,\ \pi/2]=W_2^l(-\pi/2,\ \pi/2)\times{\mathbb C}^{2m}
$$
given by
$$
\tilde{\cal L}_g(\lambda)\tilde U=\{\tilde{\cal P}(\varphi,\
D_\varphi,\ \lambda)\tilde U(\varphi),\
  \tilde{\cal B}_{i\mu 0}(\varphi,\ D_\varphi,\ \lambda)\tilde U(\varphi)|_{\varphi=(-1)^j\pi/2}\},\ j=1,\ 2.
$$

If $g\in{\cal K}_{3},$ we introduce the bounded operator
$$
 {\cal L}_g={\cal P}(D_y):H_b^{l+2m}({\mathbb R}^2)\to  H_b^l({\mathbb R}^2).
$$
Let us also introduce the bounded operator
$$
 \tilde{\cal L}_g(\lambda)=\tilde{\cal P}(\varphi,\ D_\varphi,\ \lambda):
 W_{2,2\pi}^{l+2m}(0,\ 2\pi)\to W_{2,2\pi}^l(0,\ 2\pi),
$$
where $W_{2,2\pi}^{l}(0,\ 2\pi)$ is the closure of the set of
infinitely differentiable $2\pi$-periodic functions in
$W_2^{l}(0,\ 2\pi).$

From~\cite[\S~1]{KondrTMMO67} and~\cite[\S~1]{SkMs86}, it follows
that for each $g\in{\cal K}_2\cup{\cal K}_3$ there is a
finite-meromorphic operator-valued function $\tilde{\cal
L}_g^{-1}(\lambda)$ such that {\rm (I)} its poles, maybe with the
exception of finitely many of them, belong to a double angle of
opening $<$$\pi$, containing the imaginary axis, and {\rm (II)}
for a $\lambda$ which is not a pole of $\tilde{\cal
L}_g^{-1}(\lambda)$, the operator $\tilde{\cal L}_g^{-1}(\lambda)$
is the bounded inverse for $\tilde{\cal L}_g(\lambda).$

From Theorem~1.1~\cite{KondrTMMO67} and results
of~\cite[\S~1]{SkMs86}, it follows that the operator ${\cal L}_g$
is an isomorphism if and only if the line $\Im\lambda=b+1-l-2m$
contains no poles of $\tilde{\cal L}_g^{-1}(\lambda).$

\begin{theorem}\label{thFredLinAllG}
Let Conditions~\ref{condEllipPinG}--\ref{condCoordTrans} and
\ref{condKinG} hold. Suppose that $b\in\mathbb R$ is such that for
all $g\in{\cal K}_1$ the operator $\mathcal L_g^{\mathcal G}$ is
an isomorphism and for all $g\in\mathcal K_2\cup\mathcal K_3$ the
operator $\mathcal L_g$ is an isomorphism.

Then the operator ${\bf L}:\tilde H_b^{l+2m}(G)\to \tilde {\cal
H}_b^l(G,\ \Upsilon)$ is Fredholm.
\end{theorem}
\begin{proof}
Notice that Lemmas~\ref{lAprLMain} and~\ref{lRegGK1} are true for
any $b\in\mathbb R$ for which the operators $\mathcal
L_g^{\mathcal G}$, $g\in{\cal K}_1$, are isomorphisms. Therefore,
using Lemmas~\ref{lAprLMain} and~\ref{lRegGK1}, analogously to the
proof of Theorem~3.4~\cite{SkMs86}, we can obtain an a priori
estimate~\eqref{eqAprLMain} (in the spaces $\tilde
H_b^l(\,\cdot\,)$) and construct a right regularizer.
\end{proof}

%% file: sect6.tex
\section{Index stability for nonlocal elliptic problems}\label{sectStabInd}
In this section, we study an influence of the transformations
$\omega_{is}$ upon the index of nonlocal elliptic problems. We
show that the index of the problem is determined by the linear
part of the transformations $\omega_{is}$ in a neighborhood of
${\cal K}_1$. Notice that, in the case where the support
$\bigcup\limits_{i,\,s}\omega_{is}(\bar\Upsilon_i)$ of nonlocal
terms does not intersect with the set ${\cal K}_1$ consisting of
the points of conjugation of nonlocal conditions, the index
stability for the corresponding problem was proved
in~\cite{SkJMAA91}.

{\bf 1.} Parallel to problem~(\ref{eqPinG}), (\ref{eqBinG}), we
consider the following problem:
\begin{equation}\label{eqPinGHat}
 {\bf P}(x,\ D) u=f_0(x) \quad (x\in G),
\end{equation}
\begin{equation}\label{eqBinGHat}
 \begin{array}{c}
   \hat{\bf B}_{i\mu}(x,\ D)u\equiv\sum\limits_{s=0}^{\hat S_i}(\hat B_{i\mu s}(x,\ D)u)(\hat\omega_{is}(x))|_{\Upsilon_i}=g_{i\mu}(x)\\
    (x\in \Upsilon_i;\ i=1,\ \dots,\ N_0;\ \mu=1,\ \dots,\ m).
 \end{array}
\end{equation}
Here ${\bf P}(x,\ D)$, $\hat B_{i\mu 0}(x,\ D)=B_{i\mu 0}(x,\ D)$
are the same\footnote{It suffices that only the principal
homogeneous parts of the operators ${\bf P}(x,\ D)$ and $\hat
B_{i\mu 0}(x,\ D)$ from this section and those
from~\S~\ref{sectStatement} coincide. But, for simplicity, we
assume that junior terms of the corresponding operators also
coincide.} differential operators as those
in~\S~\ref{sectStatement}, $\hat B_{i\mu s}(x,\ D)$ ($s=1,\
\dots,\ \hat S_i$) are some differential operators of orders
$m_{i\mu}$ with complex-valued $C^\infty$-coefficients;
$\hat\omega_{is}$ ($i=1,\ \dots,\ N_0;$ $s=1,\ \dots,\ \hat S_i$)
are infinitely differentiable non-degenerate transformations
mapping some neighborhood ${\cal O}_i$ of the manifold
$\Upsilon_i$ onto $\hat\omega_{is}({\cal O}_i)$ so that
$\hat\omega_{is}(\Upsilon_i)\subset G$; $\omega_{i0}(x)\equiv x$.
We assume that the set
$$
  \hat{\cal K}=
  \left\{\bigcup_i (\bar \Upsilon_i\setminus \Upsilon_i)\right\}\cup
  \left\{\bigcup_{i,\,s}
  \hat\omega_{is}(\bar \Upsilon_i\setminus \Upsilon_i)\right\}\cup
  \left\{\bigcup_{j,\,p}\,\bigcup_{i,\,s} \hat\omega_{jp}(\hat\omega_{is}(\bar \Upsilon_i\setminus \Upsilon_i)\cap\Upsilon_j)\right\}
$$
can be represented in the form $\hat{\cal
K}=\bigcup\limits_{j=1}^3\bigcup\limits_{p=1}^{\hat N_j} \hat{\cal
K}_{jp},$ where
$$
\hat{\cal K}_1=\bigcup\limits_{p=1}^{\hat N_1} \hat{\cal K}_{1p}=
\partial G\setminus\bigcup\limits_{i=1}^{N_0}
\Upsilon_i,\quad \hat{\cal K}_2=\bigcup\limits_{p=1}^{\hat N_2}
\hat{\cal K}_{2p}\subset\bigcup\limits _{i=1}^{N_0}
\Upsilon_i,\quad \hat{\cal K}_3=\bigcup\limits_{p=1}^{\hat N_3}
\hat{\cal K}_{3p}\subset G
$$
(cf.~\eqref{eqDecompK}). Here $\hat{\cal K}_{jp}$ are disjoint
$(n-2)$-dimensional $C^\infty$-manifolds without a boundary
(points if $n=2$); moreover, $\hat N_1=N_1$, $\hat{\cal
K}_{1p}={\cal K}_{1p}$, $p=1,\ \dots,\ N_1$.

Let the transformations $\hat\omega_{is}$ satisfy
Conditions~\ref{condOrb} and~\ref{condCoordTrans}. Furthermore, we
assume that the operators $\hat B_{i\mu s}(x,\ D)$ and the
transformations $\hat\omega_{is}$ ($s=1,\ \dots,\ \hat S_i$) are
such that for each $g\in\hat{\mathcal K}_1=\mathcal K_1$ the
operator $\mathcal L_g^{\hat\omega}$ (which is defined similarly
to the operator $\mathcal L_g^{\omega}$
from~\S~\ref{sectStatement}) equals the operator $\mathcal
L_g^{\mathcal G}$ defined in~\S~\ref{sectStatement}.

Thus, $\hat\omega_{is}$ is a linear part of $\omega_{is}$ in a
neighborhood of ${\cal K}_1$.

We introduce the bounded operator corresponding to nonlocal
problem~(\ref{eqPinGHat}), (\ref{eqBinGHat}):
$$
  \hat{\bf L}=\{{\bf P}(x,\ D),\ \hat{\bf B}_{i\mu}(x,\ D)\}:
 H_b^{l+2m}(G)\to  {\cal H}_b^l(G,\ \Upsilon).
$$

\begin{theorem}\label{thIndStab}
Let the conditions of Lemma~\ref{lAprLMain} hold and $b>l+2m-1$.
Then the operators ${\bf L},\ \hat{\bf L}:H_b^{l+2m}(G)\to {\cal
H}_b^l(G,\ \Upsilon)$ are Fredholm and $\ind{\bf L}=\ind\hat{\bf
L}$.
\end{theorem}
\begin{proof}
We consider the operator ${\bf L}_t: H_b^{l+2m}(G)\to {\cal
H}_b^l(G,\ \Upsilon)$ given by
$$
 {\bf L}_tu=\{{\bf P}(x,\ D)u,\ {\bf B}_{i\mu}(x,\ D)+t(\hat{\bf B}_{i\mu}(x,\ D)-{\bf B}_{i\mu}(x,\ D))\}.
$$
Obviously, ${\bf L}_0={\bf L}$, ${\bf L}_1=\hat{\bf L}$.

In a neighborhood of ${\cal K}_1$, the transformations
$\omega_{is}$ and $\hat\omega_{is}$ coincide up to infinitesimals;
therefore, by Theorem~\ref{thFredLinG}, the operators ${\bf L}_t$
are Fredholm for all $t$. Furthermore, for all $t_0$ and $t$, we
have
$$
 \| {\bf L}_tu- {\bf L}_{t_0}u\|_{{\cal H}_b^l(G,\ \Upsilon)}\le k_{t_0}|t-t_0|\cdot \|u\|_{H_b^{l+2m}(G)},
$$
where $k_{t_0}>0$ is independent of $t\in [0,\ 1]$. Hence, by
Theorem~16.2~\cite{Kr}, we have $\ind{\bf L}_t=\ind{\bf L}_{t_0}$
for all $t$ from some small neighborhood of $t_0$. These
neighborhoods cover the segment $[0,\ 1]$. Choosing a finite
subcovering, we get $\ind{\bf L}=\ind{\bf L}_0=\ind{\bf
L}_1=\ind\hat{\bf L}$.
\end{proof}

\smallskip

Analogously to the above, using Theorem~\ref{thFredLinAllG}
instead of Theorem~\ref{thFredLinG}, one can prove the index
stability for nonlocal problem~(\ref{eqPinG}), (\ref{eqBinG}) in
the case where $n=2$, $b\in\mathbb R$.

Let us suppose that $\hat N_j= N_j$, $\hat{\cal K}_{jp}={\cal
K}_{jp}$, $j=1,\ 2,\ 3$, $p=1,\ \dots,\ N_j$.
\begin{theorem}\label{thIndStab2}
Let the conditions of Theorem~\ref{thFredLinAllG} hold. Then the
operators ${\bf L},\ \hat{\bf L}:\tilde H_b^{l+2m}(G)\to \tilde
{\cal H}_b^l(G,\ \Upsilon)$ are Fredholm and $\ind{\bf
L}=\ind\hat{\bf L}$.
\end{theorem}

\medskip

{\bf 2.} In this subsection, we present another proof of
Theorem~\ref{thIndStab2}, based upon ideas of~\cite{SkJMAA91}.
(Notice that, using Lemma~5.2~\cite{KovSk}, one can similarly
prove Theorem~\ref{thIndStab}.) The proof given below is more
complicated; however it makes clear the phenomenon---{\it why
index of the operator is completely determined by the linear part
of the transformations $\omega_{is}$ in a neighborhood of ${\cal
K}_1$}. We show that if the operators ${\bf L}$ and $\hat{\bf L}$
are both Fredholm, then the restriction of their difference to the
kernel $\ker({\bf P})\subset \tilde H_b^{l+2m}(G)$ of the operator
${\bf P}={\bf P}(y,\ D)$ (we remind that $x=y$ if $n=2$) can be
``reduced'' to the sum of an operator with an arbitrary small norm
and an operator the square of which is compact. The first operator
appears at the expense of the nonlinear part of the
transformations $\omega_{is}$ near ${\cal K}_1$ while the second
one appears at the expense of transformations originating the sets
${\cal K}_2$ and ${\cal K}_3$ (see~\S~\ref{sectStatement}). Notice
that this ``reduction'' does not contradict the example
of~\S~\ref{sectEx} since the ``reduction'' procedure contains
projecting to the subspace $\ker({\bf P})$ of infinite
codimension. By the same reason, the considerations below do not
prove that the operator $\hat{\bf L}$ is Fredholm whenever ${\bf
L}$ is Fredholm (or vice versa). The only thing they imply is that
$\ind{\bf L}=\ind\hat{\bf L}$ whenever we are a priori aware of
${\bf L}$ and $\hat{\bf L}$ being both Fredholm.

Thus, let us proceed to the alternative proof of
Theorem~\ref{thIndStab2}.

1) We introduce the operators
$$
 {\bf B},\ \hat{\bf B}:\tilde H_b^{l+2m}(G)\to \tilde{\cal H}_b^l(\partial G)=
 \prod\limits_{i=1}^{N_0}\prod\limits_{\mu=1}^m \tilde H_b^{l+2m-m_{i\mu}-1/2}(\Upsilon_i)
$$
given by ${\bf B}=\{{\bf B}_{i\mu}(y,\ D)\},\ \hat{\bf
B}=\{\hat{\bf B}_{i\mu}(y,\ D)\}$. We denote by ${\bf C},\
\hat{\bf C}$ the restrictions of the operators ${\bf B},\ \hat{\bf
B}$ to the subspace $\ker({\bf P})\subset \tilde H_b^{l+2m}(G)$.
By Theorem~\ref{thFredLinG}, the operators ${\bf L},\ \hat{\bf L}$
are Fredholm. Therefore, by virtue of Lemma~1.1~\cite{SkJMAA91},
the operators ${\bf C},\ \hat{\bf C}$ are also Fredholm. Now, to
prove Theorem~\ref{thIndStab2}, it suffices to show that $\ind
{\bf C}=\ind\hat{\bf C}$.

2) We denote by ${\bf C}^1,\ \hat{\bf C}^1$ the restrictions of
${\bf C},\ \hat{\bf C}$ to the subspace $\ker({\bf
C})^\bot\subset\ker({\bf P})$. It is obvious that ${\bf C}^1={\bf
C}{\bf I}_0$, $\hat{\bf C}^1=\hat{\bf C}{\bf I}_0$, where ${\bf
I}_0:\ker({\bf C})^\bot\to\ker({\bf P})$ is the operator of
embedding of $\ker({\bf C})^\bot$ into $\ker({\bf P})$. Clearly,
we have $\dim\ker({\bf I}_0)=0$, $\codim{\cal R}({\bf
I}_0)=\dim\ker({\bf C})=m_0<\infty$. Therefore, from
Theorem~12.2~\cite{Kr}, it follows that
 \begin{gather*}
 \ind {\bf C}^1=\ind {\bf C}+\ind {\bf I}_0=\ind {\bf C}-m_0,\\
 \ind \hat{\bf C}^1=\ind \hat{\bf C}+\ind {\bf I}_0=\ind \hat{\bf C}-m_0.
 \end{gather*}
Thus, it suffices to prove that $\ind {\bf C}^1=\ind\hat{\bf
C}^1$.

3) We denote by ${\bf P}_\bot$ the operator that orthogonally
projects $\tilde {\cal H}_b^l(\partial G)$ onto ${\cal R}({\bf
C}^1)^\bot$. Since $\codim{\cal R}({\bf C}^1)<\infty$, it follows
that the operator ${\bf P}_\bot$ is finite-dimensional. Therefore,
we have
$$
 \ind\hat{\bf C}^1=\ind\big({\bf C}^1+({\bf I}-{\bf P}_\bot)(\hat{\bf C}^1-{\bf C}^1)\big).
$$
Hence, it suffices to prove that
$$
 \ind {\bf C}^1=\ind\big({\bf C}^1+({\bf I}-{\bf P}_\bot)(\hat{\bf C}^1-{\bf C}^1)\big).
$$
Since ${\bf C}^1u,\ {\bf C}^1u+({\bf I}-{\bf P}_\bot)(\hat{\bf
C}^1-{\bf C}^1)u\in{\cal R}({\bf C}^1)$ for $u\in\ker({\bf
C})^\bot$, we may regard ${\bf C}^1,\ {\bf C}^1+({\bf I}-{\bf
P}_\bot)(\hat{\bf C}^1-{\bf C}^1)$ as the operators acting from
$\ker({\bf C})^\bot$ into ${\cal R}({\bf C}^1)$. In this case, the
indices of these operators increase the same number
$m_1=\codim{\cal R}({\bf C}^1)$.

Evidently, the operator ${\bf C}^1:\ker({\bf C})^\bot\to {\cal
R}({\bf C}^1)$ has the bound inverse ${\bf R}_1=({\bf
C}^1)^{-1}:{\cal R}({\bf C}^1)\to \ker({\bf C})^\bot$ and
$\ind{\bf C}^1=0$. By Theorem~12.2 \cite{Kr}, we have
$$
\ind\big({\bf C}^1+({\bf I}-{\bf P}_\bot)(\hat{\bf C}^1-{\bf
C}^1)\big)=\ind\big({\bf I}+{\bf R}_1({\bf I}-{\bf
P}_\bot)(\hat{\bf C}^1-{\bf C}^1)\big).
$$
It remains to show that $\ind\big({\bf I}+{\bf R}_1({\bf I}-{\bf
P}_\bot)(\hat{\bf C}^1-{\bf C}^1)\big)=0$.

4) Let us introduce a function $\psi_{\varepsilon}\in
C_0^\infty({\mathbb R}^2)$ such that $\psi_{\varepsilon}(y)=1$ for
$y\in{\cal O}_{\varepsilon/2}({\cal K})$,
$\psi_{\varepsilon}(y)=0$ for $y\notin{\cal O}_{\varepsilon}({\cal
K})$, and
\begin{equation}\label{eqPsiinG}
|D^\alpha\psi_{\varepsilon}(y)|\le
k_{\alpha}(\tilde\rho(y))^{-|\alpha|}\quad \big(y\in {\cal
O}_{\varepsilon}({\cal K})\big),
\end{equation}
where $k_{\alpha}>0$ is independent of $\varepsilon$.

We consider the operators ${\bf A}_1,\ {\bf A}_2:\ker({\bf
C})^\bot\to\ker({\bf C})^\bot$ given by
 \begin{gather*}
 {\bf A}_1u={\bf R}_1({\bf I}-{\bf P}_\bot)(\hat{\bf B}-{\bf B})\psi_\varepsilon u,\\
 {\bf A}_2u={\bf R}_1({\bf I}-{\bf P}_\bot)(\hat{\bf B}-{\bf B})(1-\psi_\varepsilon)u.
 \end{gather*}
It is clear that ${\bf I}+ {\bf A}_1+ {\bf A}_2={\bf I}+{\bf
R}_1({\bf I}-{\bf P}_\bot)(\hat{\bf C}^1-{\bf C}^1)$. Since the
support of $(1-\psi_\varepsilon)u$ does not intersect with the
origin, it follows from the proof of Theorem~3.1~\cite{SkJMAA91}
that the operator $({\bf A}_2)^2$ is compact.

Let us study the operator ${\bf A}_1$. Since the operator ${\bf
R}_1({\bf I}-{\bf P}_\bot)$ is bounded, it follows that
$$
 \|{\bf A}_1u\|_{\tilde H_b^{l+2m}(G)}\le c\|(\hat{\bf B}-{\bf B})\psi_\varepsilon u\|_{\tilde {\cal H}_b^l(\partial G)}.
$$
From this, using the unity partition method and
estimates~(\ref{eqAprL2})--(\ref{eqAprL7}), followed
by~(\ref{eqPsiinG}), we obtain
\begin{multline}\label{eqIndStab'1}
 \|{\bf A}_1u\|_{\tilde H_b^{l+2m}(G)}\le c_1(\varepsilon\|\psi_\varepsilon u\|_{\tilde H_b^{l+2m}(G)}+\|{\bf P}\psi_\varepsilon u\|_{\tilde H_b^l(G)})+
 k_1(\varepsilon)\|u\|_{\tilde H_b^{l+2m-1}(G)}\le\\
\le c_2(\varepsilon\|u\|_{\tilde H_b^{l+2m}(G)}+\|{\bf
P}\psi_\varepsilon u\|_{\tilde H_b^l(G)})+
  k_1(\varepsilon)\|u\|_{\tilde H_b^{l+2m-1}(G)}.
\end{multline}
Since $u\in\ker({\bf P})$, from~(\ref{eqIndStab'1}) and Leibniz'
formula, we get
\begin{equation}\label{eqIndStab'2}
 \|{\bf A}_1u\|_{\tilde H_b^{l+2m}(G)}\le c_2\varepsilon\|u\|_{\tilde H_b^{l+2m}(G)}+k_2(\varepsilon)\|u\|_{\tilde H_b^{l+2m-1}(G)},
\end{equation}
where $c_2$ is independent of $\varepsilon$.
From~(\ref{eqIndStab'2}), the compactness of the embedding $\tilde
H_b^{l+2m}(G)\subset \tilde H_b^{l+2m-1}(G)$, and
Lemma~\ref{lSmallComp}, it follows that ${\bf A}_1={\bf M}_1+{\bf
F}_1$, where $\|{\bf M}_1\|\le 2c_2\varepsilon$ and the operator
${\bf F}_1$ is finite-dimensional.

Thus, we have ${\bf R}_1({\bf I}-{\bf P}_\bot)(\hat{\bf C}^1-{\bf
C}^1)={\bf M}_1+{\bf F}_1+{\bf A}_2$. Therefore, choosing
sufficiently small $\varepsilon$, we obtain from Theorems~15.4
and~16.2~\cite{Kr} that $\ind\big({\bf I}+{\bf R}_1({\bf I}-{\bf
P}_\bot)(\hat{\bf C}^1-{\bf C}^1)\big)=0$.

\bigskip

The author is very grateful to A.~L.~Skubachevskii for the
statement of the problem and attention to this work.